\documentclass[10pt]{amsart}
\usepackage[cp1251]{inputenc}
\usepackage[english]{babel}
\usepackage{amsmath}
\usepackage{amssymb}
\usepackage{amsfonts}

\begin{document}
\renewcommand{\refname}{References}

\thispagestyle{empty}

\title[Numerical Simulation of 2.5-Set of Iterated 
Stratonovich Stochastic Integrals]
{Numerical Simulation of 2.5-Set of Iterated Stratonovich 
Stochastic Integrals of Multiplicities 1 to 5 From the 
Taylor--Stratonovich Expansion}
\author[D.F. Kuznetsov]{Dmitriy F. Kuznetsov}
\address{Dmitriy Feliksovich Kuznetsov
\newline\hphantom{iii} Peter the Great Saint-Petersburg Polytechnic University,
\newline\hphantom{iii} Polytechnicheskaya ul., 29,
\newline\hphantom{iii} 195251, Saint-Petersburg, Russia}%
\email{sde\_kuznetsov@inbox.ru}

\vspace{5mm}

\maketitle{\small
\begin{quote}
\noindent{\sc Abstract.} 
The article is devoted to construction of effective procedures 
of the mean-square approximation for iterated Stratonovich stochastic 
integrals of multiplicities 1 to 5. We apply the method of 
generalized multiple Fourier series for approximation of 
iterated stochastic integrals. More precisely, we use multiple 
Fourier--Legendre series converging in the sense of norm
in Hilbert space $L_2([t,T]^k),$ $k\in\mathbb{N}.$
Considered iterated Stra\-to\-no\-vich stochastic integrals are part 
of the Taylor--Stratonovich expansion. That is why the results of 
the article can be applied to implementation of 
numerical methods with the orders 1.0, 1.5, 2.0 and 2.5
of strong convergence 
for Ito stochastic differential 
equations with multidimensional non-commutative noise.
\medskip
\end{quote}
}

\vspace{3mm}


\setlength{\baselineskip}{2.0em}

\tableofcontents

\setlength{\baselineskip}{1.2em}


\section{Introduction}

\vspace{5mm}

Let $(\Omega,$ ${\rm F},$ ${\sf P})$ be a complete probability space, let 
$\{{\rm F}_t, t\in[0,T]\}$ be a nondecreasing 
right-continous family of $\sigma$-algebras of ${\rm F},$
and let ${\bf f}_t$ be a standard $m$-dimensional Wiener 
stochastic process, which is
${\rm F}_t$-measurable for any $t\in[0, T].$ We assume that the components
${\bf f}_{t}^{(i)}$ $(i=1,\ldots,m)$ of this process are independent. Consider
an Ito stochastic differential equation (SDE) in the integral form

\vspace{-1mm}
\begin{equation}
\label{1.5.2}
{\bf x}_t={\bf x}_0+\int\limits_0^t {\bf a}({\bf x}_{\tau},\tau)d\tau+
\int\limits_0^t B({\bf x}_{\tau},\tau)d{\bf f}_{\tau},\ \ \
{\bf x}_0={\bf x}(0,\omega).
\end{equation}

\vspace{2mm}
\noindent
Here ${\bf x}_t$ is some $n$-dimensional stochastic process 
satisfying to the equation (\ref{1.5.2}). 
The nonrandom functions ${\bf a}: \mathbb{R}^n\times[0, T]\to\mathbb{R}^n$,
$B: \mathbb{R}^n\times[0, T]\to\mathbb{R}^{n\times m}$
guarantee the existence and uniqueness up to stochastic equivalence 
of a solution
of the equation (\ref{1.5.2}) \cite{1}. The second integral on the 
right-hand side of (\ref{1.5.2}) is 
interpreted as an Ito stochastic integral.
Let ${\bf x}_0$ be an $n$-dimensional random variable, which is 
${\rm F}_0$-measurable and 
${\sf M}\{\left|{\bf x}_0\right|^2\}<\infty$ 
(${\sf M}$ denotes a mathematical expectation).
We assume that
${\bf x}_0$ and ${\bf f}_t-{\bf f}_0$ are independent when $t>0.$

It is well known \cite{KlPl2}-\cite{KPS}
that Ito SDEs are 
adequate mathematical models of dynamic systems under 
the influence of random disturbances. One of the effective approaches 
to numerical integration of 
Ito SDEs is an approach based on 
the Taylor--Ito and 
Taylor--Stratonovich expansions
\cite{KlPl2}-\cite{xxx333}. 
The most important feature of such 
expansions is a presence in them of the so-called iterated
Ito and Stratonovich stochastic integrals, which play the key 
role for solving the 
problem of numerical integration of Ito SDEs
and have the 
following form

\vspace{-1mm}
\begin{equation}
\label{ito}
J[\psi^{(k)}]_{T,t}=\int\limits_t^T\psi_k(t_k) \ldots \int\limits_t^{t_{2}}
\psi_1(t_1) d{\bf w}_{t_1}^{(i_1)}\ldots
d{\bf w}_{t_k}^{(i_k)},
\end{equation}

\begin{equation}
\label{str}
J^{*}[\psi^{(k)}]_{T,t}=
{\int\limits_t^{*}}^T
\psi_k(t_k) \ldots 
{\int\limits_t^{*}}^{t_2}
\psi_1(t_1) d{\bf w}_{t_1}^{(i_1)}\ldots
d{\bf w}_{t_k}^{(i_k)},
\end{equation}

\vspace{3mm}
\noindent
where $\psi_1(\tau),\ldots,\psi_k(\tau)$ are continuous nonrandom
functions 
on $[t,T],$ ${\bf w}_{\tau}^{(i)}={\bf f}_{\tau}^{(i)}$
for $i=1,\ldots,m$ and
${\bf w}_{\tau}^{(0)}=\tau,$\ \
$i_1,\ldots,i_k = 0, 1,\ldots,m,$

\vspace{-1mm}
$$
\int\limits\ \hbox{and}\ \int\limits^{*}
$$ 

\vspace{3mm}
\noindent
denote Ito and 
Stratonovich stochastic integrals,
respectively (in this paper, 
we use the definition of the Stratonovich stochastic integral from \cite{KlPl2}).

Note that $\psi_l(\tau)\equiv 1$ $(l=1,\ldots,k)$ and
$i_1,\ldots,i_k = 0, 1,\ldots,m$ in the classical Taylor--Ito
and Taylor--Stratonovich expansions
\cite{KlPl2}-\cite{KlPl1}. At the same time
$\psi_l(\tau)\equiv (t-\tau)^{q_l}$ ($l=1,\ldots,k$; 
$q_1,\ldots,q_k=0, 1, 2,\ldots $) and  $i_1,\ldots,i_k = 1,\ldots,m$ 
in the unified Taylor--Ito
and Taylor--Stratonovich expansions
\cite{kk5}-\cite{xxx333}.

Effective solution 
of the problem of
combined mean-square approximation of collections 
of the iterated Ito and Stratonovich stochastic integrals
(\ref{ito}), (\ref{str}) of multiplicities 1 to 5 and beyond
composes the subject of the article.

We want to mention in short that there are 
two main criteria of numerical methods convergence 
for Ito SDEs \cite{KlPl2}-\cite{Mi3}:  
a strong or mean-square
criterion and a 
weak criterion where the subject of approximation is not the solution 
of Ito SDE, simply stated, but the 
distribution of Ito SDE solution.

Using the strong numerical methods, we can build
sample pathes
of Ito SDEs numerically. 
These methods require the combined mean-square approximation of collections 
of the iterated Ito and Stratonovich stochastic integrals
(\ref{ito}) and (\ref{str}). 

The strong numerical methods are using when constructing new mathematical 
models on the basis of Ito SDEs, when
solving the filtering problem of signal under the influence
of random disturbance in various arrangements, 
when solving the problem of stochastic 
optimal control, when solving the problem
of testing procedures of evaluating parameters of stochastic 
systems etc. \cite{KlPl2}-\cite{KPS}.

The problem of effective jointly numerical modeling 
(in accordance to the mean-square convergence criterion) of the iterated 
Ito and Stratonovich stochastic integrals 
(\ref{ito}) and (\ref{str}) is 
difficult from 
theoretical and computing point of view \cite{KlPl2}-\cite{KPS},
\cite{2006}-\cite{Rybakov1000}.

The only exception is connected with the narrow particular case, when 
$i_1=\ldots=i_k\ne 0$ and
$\psi_1(s),\ldots,\psi_k(s)\equiv \psi(s)$.
This case allows 
the investigation with using of the Ito formula 
\cite{KlPl2}-\cite{Mi3}.

Note that even for the mentioned coincidence ($i_1=\ldots=i_k\ne 0$), 
but for different 
functions $\psi_1(s),\ldots,\psi_k(s)$ the mentioned 
difficulties persist, and 
relatively simple families of 
iterated Ito and Stratonovich stochastic integrals, 
which can be often 
met in the applications, cannot be represented effectively in a finite 
form (for the mean-square approximation) using the system of standard 
Gaussian random variables.

Note that for a number of special types of Ito SDEs 
the problem of approximation of iterated
stochastic integrals can be simplified but cannot be solved. The equations
with additive vector noise, with additive scalar or non-additive scalar 
noise, with a small parameter are related to such 
types of equations \cite{KlPl2}-\cite{Mi3}. 
For the mentioned types of equations, simplifications 
are connected with the fact that either some coefficient functions 
from stochastic analogues of the Taylor formula 
(Taylor--Ito and Taylor--Stratonovich expansions)
identically equal to zero, 
or scalar noise has an essential effect, or due to the presence 
of a small parameter we can neglect some members from stochastic 
analogues of the Taylor formula, which include difficult for approximation 
iterated stochastic integrals \cite{KlPl2}-\cite{Mi3}.
In this article, we consider Ito SDEs 
with multidimentional and non-additive noise. 
The conditions of commutativity of the noise \cite{KlPl2}
are also not used.

Seems that iterated stochastic integrals can be approximated by multiple 
integral sums of different types \cite{Mi2}, \cite{Mi3}, \cite{allen}. 
However, this approach implies partitioning of the interval 
of integration $[t, T]$ of iterated stochastic integrals 
(the length $T-t$ of this interval is a small 
value, because it is a step of integration of numerical methods for 
Ito SDEs) and according to numerical 
experiments this additional partitioning leads to significant calculating 
costs \cite{2006}.

In \cite{Mi2} (also see \cite{KlPl2}, \cite{Mi3})
Milstein G.N. proposed to expand (\ref{ito}) or (\ref{str})
into iterated series of products
of standard Gaussian random variables by representing the Wiener
process as a trigonometric Fourier series with random coefficients 
(the version of the so-called Karhunen--Loeve expansion for
the Brownian bridge process).
For example, 
to obtain the Milstein expansion of (\ref{str}), the truncated Fourier
expansions of components of the Wiener process ${\bf f}_s$ must be
iteratively substituted in the single integrals, and the integrals
must be calculated, starting from the innermost integral.
This is a complicated procedure that does not lead to a general
expansion of (\ref{str}) valid for an arbitrary multiplicity $k.$
For this reason, only expansions of single, double, and triple
stochastic integrals (\ref{ito}) and (\ref{str}) were presented 
in \cite{KlPl2} (the integrals (\ref{str}) for $k=1, 2, 3$)
and in \cite{Mi2}, \cite{Mi3} (the integrals (\ref{ito}) for $k=1, 2$) 
for the simplest case $\psi_1(s), \psi_2(s), \psi_3(s)\equiv 1;$ 
$i_1, i_2, i_3=0,1,\ldots,m.$ Moreover, the Milstein 
approach \cite{Mi2} leads to iterated application
of the operation of limit transition (see above).

It should be noted that the authors of the works
\cite{KlPl2}
(Sect.~5.8, pp.~202--204), \cite{KPS} (pp.~82-84),
\cite{KPW} (pp.~438-439),  
\cite{Zapad-9} (pp.~263-264) use 
the Wong--Zakai approximation 
\cite{W-Z-1}-\cite{Watanabe} (without rigorous proof) within the frames
of the method of expansion of iterated stochastic integrals
\cite{Mi2} (1988) based on the series expansion 
of the Brownian bridge process (version
of the so-called Karhunen-Loeve expansion).
See discussions in \cite{2018a} (Sect.~2.18, 6.2), 
\cite{xxx333} (Sect.~2.6.2, 6.2)
\cite{arxiv-1} (Sect.~11),
\cite{arxiv-3} (Sect.~8),
\cite{arxiv-4} (Sect.~11),
\cite{arxiv-5} (Sect.~6),
\cite{arxiv-7} (Sect.~6)
for detail.

Note that in \cite{rr} the method of expansion
of iterated (double) Ito stochastic integrals (\ref{ito}) 
($k=2;$ $\psi_1(s), \psi_2(s) \equiv 1;$ $i_1, i_2 =1,\ldots,m$) 
based on expansion
of the Wiener process using Haar functions and 
trigonometric functions has been considered.
The restrictions of the method \cite{rr} are also connected
with iterated application of the operation
of limit transition (as in the Milstein approach \cite{Mi2} (1988)) 
at least starting from the third 
multiplicity of iterated stochastic integrals.

It is necessary to note that the Milstein approach \cite{Mi2} 
excelled in several times or even in several orders
the methods based on multiple integral sums 
\cite{Mi2}, \cite{Mi3}, \cite{allen}
considering computational costs in the sense 
of their diminishing.

An alternative strong approximation method was 
proposed for (\ref{str}) in \cite{300a}, \cite{400a} 
(also see \cite{2011-2}-\cite{xxx333},
\cite{2010-1}-\cite{2013}, \cite{arxiv-23}),
where $J^{*}[\psi^{(k)}]_{T,t}$ was represented as the multiple stochastic 
integral
from the certain discontinuous nonrandom function of $k$ variables, and the 
function
was then expressed as the iterated generalized Fourier series in complete
systems of continuous functions that are orthonormal in the space
$L_2([t, T]).$
In \cite{300a}, \cite{400a} 
(also see \cite{2011-2}-\cite{xxx333},
\cite{2010-1}-\cite{2013}, \cite{arxiv-23}) 
the cases of Legendre polynomials and
trigonometric functions are considered in detail.
As a result,
the general iterated series expansion of (\ref{str}) in terms of products
of standard Gaussian random variables was obtained in 
\cite{300a}, \cite{400a} 
(also see \cite{2011-2}-\cite{xxx333},
\cite{2010-1}-\cite{2013}, \cite{arxiv-23})
for an arbitrary multiplicity $k.$
Hereinafter, this method is referred to as the method of generalized
iterated Fourier series.

It was shown in \cite{300a}, \cite{400a} 
(also see \cite{2011-2}-\cite{xxx333},
\cite{2010-1}-\cite{2013}, \cite{arxiv-23}) that the method of 
generalized iterated Fourier series leads to the 
Milstein expansion \cite{Mi2} of (\ref{str})
in the case of trigonometric 
functions
and to a substantially simpler expansion of (\ref{str}) in the case
of Legendre polynomials.

Note that the method of generalized 
iterated Fourier series as well as the Milstein approach
\cite{Mi2}
lead to iterated application of the operation of limit transition. 
As mentioned above, this problem appears for iterated (triple) 
stochastic integrals
($i_1, i_2, i_3=1,\ldots,m$)
or even for some iterated (double) stochastic integrals 
in the case, when $\psi_1(s),$ 
$\psi_2(s)\not\equiv 1$ ($i_1, i_2=1,\ldots,m$)
\cite{2006} (also see \cite{2011-2}-\cite{200a},
\cite{301a}-\cite{arxiv-12},
\cite{arxiv-24}-\cite{arxiv-6}).
The mentioned problem (iterated application of the operation 
of limit transition) not appears 
in the efficient method, which 
is considered for (\ref{ito}) in Theorems 1, 2 (see below)
\cite{2006}-\cite{200a},
\cite{301a}-\cite{arxiv-12},
\cite{arxiv-24}-\cite{new-art-1-xxy}.

The idea of this method is as follows: 
the iterated Ito stochastic 
integral (\ref{ito}) of multiplicity $k$ is represented as 
the multiple stochastic 
integral from the certain discontinuous nonrandom function of $k$ variables
defined on the hypercube $[t, T]^k$, where $[t, T]$ is the interval of 
integration of the iterated Ito stochastic 
integral (\ref{ito}). Then, 
the indicated 
nonrandom function is expanded in the hypercube $[t, T]^k$
into the generalized 
multiple Fourier series converging 
in the mean-square sense
in the space 
$L_2([t,T]^k)$. After a number of nontrivial transformations we come 
(see Theorems 1, 2 below) to the 
mean-square convergening expansion of 
the iterated Ito stochastic 
integral (\ref{ito})
into the multiple 
series of products
of standard  Gaussian random 
variables. The coefficients of this 
series are the coefficients of 
generalized multiple Fourier series for the mentioned nonrandom function 
of $k$ variables, which can be calculated using the explicit formula 
regardless of the multiplicity $k$ of the
iterated Ito stochastic 
integral (\ref{ito}).
Hereinafter, this method is referred to as the method of generalized
multiple Fourier series.

Thus, we obtain the following useful possibilities
of the method of generalized multiple Fourier series.

\vspace{2mm}

1. There is an explicit formula (see (\ref{ppppa}) below) for calculation 
of expansion coefficients 
of the iterated Ito stochastic integral (\ref{ito}) with any
fixed multiplicity $k$. 

\vspace{2mm}

2. We have possibilities for exact calculation of the mean-square 
error of approximation 
of the iterated Ito stochastic integral (\ref{ito})
\cite{2017}-\cite{xxx333}, \cite{17a}, \cite{arxiv-2}.

\vspace{2mm}

3. Since the used
multiple Fourier series is a generalized in the sense
that it is constructed using various complete orthonormal
systems of functions in the space $L_2([t, T])$, then we 
have new possibilities 
for approximation --- we can 
use not only trigonometric functions as in \cite{KlPl2}-\cite{Mi3},
but Legendre polynomials.

\vspace{2mm}

4. As it turned out \cite{2006}-\cite{200a},
\cite{301a}-\cite{arxiv-12},
\cite{arxiv-24}-\cite{arxiv-6xx} it is more convenient to work 
with Legendre polynomials for constructing of approximations 
of the iterated Ito stochastic integrals (\ref{ito}). 
Approximations based on the Legendre polynomials are essentially simpler 
than their analogues based on the trigonometric functions.
Another advantages of the application of Legendre polynomials 
in the framework of the mentioned problem are considered
in \cite{2018a}-\cite{xxx333}, \cite{29a}, \cite{301a}.
 
\vspace{2mm}

5. An approach based on the Karhunen--Loeve expansion
of the Brownian bridge process (also see \cite{rr})
leads to 
iterated application of the operation of limit
transition (the operation of limit transition 
is implemented only once in Theorems 1, 2 (see below))
starting from  
the second multiplicity (in the general case) 
and third multiplicity (for the case
$\psi_1(s), \psi_2(s), \psi_3(s)\equiv 1;$ 
$i_1, i_2, i_3=1,\ldots,m$)
of iterated Ito stochastic integrals.
Multiple series (the operation of limit transition 
is implemented only once) are more convenient 
for approximation than the iterated ones
(iterated application of the operation of limit
transition), 
since partial sums of multiple series converge for any possible case of  
convergence to infinity of their upper limits of summation 
(let us denote them as $p_1,\ldots, p_k$). 
For example,
when $p_1=\ldots=p_k=p\to\infty$. 
For iterated series, the condition $p_1=\ldots=p_k=p\to\infty$ obviously 
does not guarantee the convergence of this series.
However, the authors of the works
\cite{KlPl2}
(Sect.~5.8, pp.~202--204), \cite{KPS} (pp.~82-84),
\cite{KPW} (pp.~438-439),  
\cite{Zapad-9} (pp.~263-264) use 
the condition $p_1=p_2=p_3=p\to\infty$
together with the Wong--Zakai approximation 
\cite{W-Z-1}-\cite{Watanabe} (but without rigorous proof) within the frames
of the method of expansion of iterated stochastic integrals
\cite{Mi2} (1988) based on the series expansion 
of the Brownian bridge process.
See discussions in \cite{2018a} (Sect.~2.18, 6.2), 
\cite{xxx333} (Sect.~2.6.2, 6.2),
\cite{arxiv-1} (Sect.~11),
\cite{arxiv-3} (Sect.~8),
\cite{arxiv-4} (Sect.~11),
\cite{arxiv-5} (Sect.~6),
\cite{arxiv-7} (Sect.~6)
for detail.

\vspace{2mm}

As it turned out, Theorems 1, 2 can be adapted for the iterated
Stratonovich stochastic integrals (\ref{str}) at least
for multiplicities 1 to 6 \cite{2011-2}-\cite{xxx333}, 
\cite{2010-2}-\cite{2013}, \cite{30a}, \cite{300a},
\cite{400a}, \cite{271a},
\cite{arxiv-4}-\cite{arxiv-8}, \cite{arxiv-23}, \cite{arxiv-6}, \cite{new-art-1-xxy}.
Expansions of these iterated Stratonovich 
stochastic integrals turned out
much simpler (see Theorems 4--10 below), than the appropriate expansions
of the iterated Ito stochastic integrals (\ref{ito}) from Theorems 1, 2.

\vspace{5mm}

\section{Explicit One-Step Strong Numerical Schemes With Orders 2.0 and 2.5 for Ito SDEs
Based on the Unified Taylor--Stratonovich expansion}

\vspace{5mm}

Consider the partition $\{\tau_j\}_{j=0}^N$ of the interval $[0, T]$ such that

\vspace{-1mm}
$$
0=\tau_0<\ldots <\tau_N=T,\ \ \
\Delta_N=
\hbox{\vtop{\offinterlineskip\halign{
\hfil#\hfil\cr
{\rm max}\cr
$\stackrel{}{{}_{0\le j\le N-1}}$\cr
}} }\Delta\tau_j,\ \ \ \Delta\tau_j=\tau_{j+1}-\tau_j.
$$

\vspace{3mm}

Let ${\bf y}_{\tau_j}\stackrel{\sf def}{=}
{\bf y}_{j},$\ $j=0, 1,\ldots,N$ be a time discrete approximation
of the process ${\bf x}_t,$ $t\in[0,T],$ which is a solution of the Ito
SDE (\ref{1.5.2}). 

\vspace{2mm}

{\bf Definiton 1}\ \cite{KlPl2}.\
{\it We will say that a time discrete approximation 
${\bf y}_{j}$\ $(j=0, 1,\ldots,N)$
corresponding to the maximal step of discretization $\Delta_N,$
converges strongly with order
$\gamma>0$ at time moment 
$T$ to the process ${\bf x}_t,$ $t\in[0,T]$,
if there exists a constant $C>0,$ which does not depend on 
$\Delta_N,$ and a $\delta>0$ such that 

\vspace{-1mm}
$$
{\sf M}\{|{\bf x}_T-{\bf y}_T|\}\le
C(\Delta_N)^{\gamma}
$$

\vspace{3mm}
\noindent
for each $\Delta_N\in(0, \delta).$}

\vspace{2mm}

Consider the explicit one-step strong numerical scheme with order 2.5 for Ito SDEs
based on the so-called unified Taylor--Stratonovich 
expansion \cite{kk6}-\cite{2010-1}, \cite{arxiv-6xx}

\vspace{2mm}
$$
{\bf y}_{p+1}={\bf y}_p+\sum_{i_{1}=1}^{m}B_{i_{1}}
\hat I_{(0)\tau_{p+1},\tau_p}^{*(i_{1})}+\Delta \bar {\bf a}
+\sum_{i_{1},i_{2}=1}^{m}G_{i_{2}}
B_{i_{1}}\hat  I_{(00)\tau_{p+1},\tau_p}^{*(i_{2}i_{1})}+
$$

\vspace{1mm}
$$
+
\sum_{i_{1}=1}^{m}\Biggl(G_{i_{1}}\bar {\bf a}\left(
\Delta \hat I_{(0)\tau_{p+1},\tau_p}^{*(i_{1})}+
\hat I_{(1)\tau_{p+1},\tau_p}^{*(i_{1})}\right)
-\bar LB_{i_{1}}\hat I_{(1)\tau_{p+1},\tau_p}^{*(i_{1})}\Biggr)+
$$

\vspace{1mm}
$$
+\sum_{i_{1},i_{2},i_{3}=1}^{m} G_{i_{3}}G_{i_{2}}
B_{i_{1}}\hat I_{(000)\tau_{p+1},\tau_p}^{*(i_{3}i_{2}i_{1})}+
\frac{\Delta^2}{2}\bar L\bar {\bf a}+
$$

\vspace{1mm}
$$
+\sum_{i_{1},i_{2}=1}^{m}
\Biggl(G_{i_{2}}\bar LB_{i_{1}}\left(
\hat I_{(10)\tau_{p+1},\tau_p}^{*(i_{2}i_{1})}-
\hat I_{(01)\tau_{p+1},\tau_p}^{*(i_{2}i_{1})}
\right)
-\bar LG_{i_{2}}B_{i_{1}}\hat I_{(10)\tau_{p+1},\tau_p}^{*(i_{2}i_{1})}
+\Biggr.
$$

\vspace{1mm}
$$
\Biggl.+G_{i_{2}}G_{i_{1}}\bar {\bf a}\left(
\hat I_{(01)\tau_{p+1},\tau_p}
^{*(i_{2}i_{1})}+\Delta \hat I_{(00)\tau_{p+1},\tau_p}^{*(i_{2}i_{1})}
\right)\Biggr)+
$$

\vspace{1mm}
$$
+
\sum_{i_{1},i_{2},i_{3},i_{4}=1}^{m}G_{i_{4}}G_{i_{3}}G_{i_{2}}
B_{i_{1}}\hat I_{(0000)\tau_{p+1},\tau_p}^{*(i_{4}i_{3}i_{2}i_{1})}+
\frac{\Delta^3}{6}LL{\bf a}+
$$

\vspace{1mm}
$$
+\sum_{i_{1}=1}^{m}\Biggl(G_{i_{1}}\bar L\bar {\bf a}\left(\frac{1}{2}
\hat I_{(2)\tau_{p+1},\tau_p}^{*(i_{1})}+
\Delta \hat I_{(1)\tau_{p+1},\tau_p}^{*(i_{1})}+
\frac{\Delta^2}{2}\hat I_{(0)\tau_{p+1},\tau_p}^{*(i_{1})}\right)\Biggr.+
$$

\vspace{1mm}
$$
\Biggl.+\frac{1}{2}\bar L 
\bar L B_{i_{1}}\hat I_{(2)\tau_{p+1},\tau_p}^{*(i_{1})}-
LG_{i_{1}}\bar {\bf a}\left(\hat I_{(2)\tau_{p+1},\tau_p}^{*(i_{1})}+
\Delta \hat I_{(1)\tau_{p+1},\tau_p}^{*(i_{1})}\right)\Biggr)+
$$

\vspace{1mm}
$$
+
\sum_{i_{1},i_{2},i_{3}=1}^m\Biggl(
G_{i_{3}}\bar LG_{i_{2}}B_{i_{1}}
\left(\hat I_{(100)\tau_{p+1},\tau_p}
^{*(i_{3}i_{2}i_{1})}-\hat I_{(010)\tau_{p+1},\tau_p}
^{*(i_{3}i_{2}i_{1})}\right)
\Biggr.+
$$

\vspace{1mm}
$$
+G_{i_{3}}G_{i_{2}}\bar LB_{i_{1}}\left(
\hat I_{(010)\tau_{p+1},\tau_p}^{*(i_{3}i_{2}i_{1})}-
\hat I_{(001)\tau_{p+1},\tau_p}^{*(i_{3}i_{2}i_{1})}\right)+
$$

\vspace{1mm}

$$
+
G_{i_{3}}G_{i_{2}}G_{i_{1}}\bar {\bf a}
\left(\Delta \hat I_{(000)\tau_{p+1},\tau_p}^{*(i_{3}i_{2}i_{1})}+
\hat I_{(001)\tau_{p+1},\tau_p}^{*(i_{3}i_{2}i_{1})}\right)
-
$$

\vspace{1mm}
$$
\Biggl.-\bar LG_{i_{3}}G_{i_{2}}B_{i_{1}}
\hat I_{(100)\tau_{p+1},\tau_p}^{*(i_{3}i_{2}i_{1})}\Biggr)+
$$

\vspace{1mm}
\begin{equation}
\label{4.470}
+\sum_{i_{1},i_{2},i_{3},i_{4},i_{5}=1}^m
G_{i_{5}}G_{i_{4}}G_{i_{3}}G_{i_{2}}B_{i_{1}}
\hat I_{(00000)\tau_{p+1},\tau_p}^{*(i_{5}i_{4}i_{3}i_{2}i_{1})},
\end{equation}

\vspace{5mm}
\noindent
where $\Delta=T/N$ $(N>1)$ is a constant (for simplicity)
step of integration,\
$\tau_p=p\Delta$ $(p=0, 1,\ldots,N)$,\
$\hat I_{(l_1\ldots l_k)s,t}^{*(i_1\ldots i_k)}$ is an
approximation of the iterated
Stratonovich stochastic integral

\vspace{-1mm}
\begin{equation}
\label{str11}
I_{(l_1\ldots \hspace{0.2mm}l_k)s,t}^{*(i_1\ldots i_k)}
=
{\int\limits_t^{*}}^s
(t-t_k)^{l_k} \ldots {\int\limits_t^{*}}^{t_{2}}
(t-t_1)^{l_1} d{\bf f}_{t_1}^{(i_1)}\ldots
d{\bf f}_{t_k}^{(i_k)},
\end{equation}

\vspace{3mm}
\noindent
where $i_1,\ldots, i_k=1,\dots,m,$\ \  $l_1,\ldots,l_k=0, 1, 2,$\ \
$k=1, 2,\ldots, 5,$

\vspace{3mm}
$$
\bar{\bf a}({\bf x},t)={\bf a}({\bf x},t)-
\frac{1}{2}\sum\limits_{j=1}^m G_jB_j({\bf x},t),
$$

\vspace{3mm}
$$
\bar L=L-\frac{1}{2}\sum\limits_{j=1}^m G_j G_j,
$$

\vspace{3mm}
$$
L= {\partial \over \partial t}
+ \sum^ {n} _ {i=1} {\bf a}_i ({\bf x},  t) 
{\partial  \over  \partial  {\bf  x}_i}
+ {1\over 2} \sum^ {m} _ {j=1} \sum^ {n} _ {l,i=1}
B_{lj} ({\bf x}, t) B_{ij} ({\bf x}, t) {\partial
^{2} \over \partial {\bf x}_l \partial {\bf x}_i},
$$

\vspace{3mm}
$$
G_i = \sum^ {n} _ {j=1} B_{ji} ({\bf x}, t)
{\partial  \over \partial {\bf x}_j}\ ,\ \ \ i=1,\ldots,m,
$$

\vspace{5mm}
\noindent
$B_i$ and $B_{ij}$ are the $i$th column and the $ij$th
element of the matrix function $B$,
${\bf a}_i$ is the $i$th element of the vector function ${\bf a},$
${\bf x}_i$ is the $i$th element
of the column ${\bf x}$, 
the functions  

\vspace{-1mm}
$$
B_{i_{1}},\  \bar {\bf a},\ G_{i_{2}}B_{i_{1}},\
G_{i_{1}}\bar {\bf a},\ \bar LB_{i_{1}},\ G_{i_{3}}G_{i_{2}}B_{i_{1}},\ 
\bar L\bar {\bf a},\ LL{\bf a},\
G_{i_{2}}\bar LB_{i_{1}},\ 
$$

\vspace{-3mm}
$$
\bar LG_{i_{2}}B_{i_{1}},\ G_{i_{2}}G_{i_{1}}\bar{\bf a},\
G_{i_{4}}G_{i_{3}}G_{i_{2}}B_{i_{1}},\ G_{i_{1}}\bar L\bar{\bf a},\
\bar L\bar LB_{i_{1}},\ \bar LG_{i_{1}}\bar {\bf a},\ 
G_{i_{3}}\bar LG_{i_{2}}B_{i_{1}},\
G_{i_{3}}G_{i_{2}}\bar LB_{i_{1}},\
$$

\vspace{-3mm}
$$
G_{i_{3}}G_{i_{2}}G_{i_{1}}\bar {\bf a},\
\bar LG_{i_{3}}G_{i_{2}}B_{i_{1}},\ 
G_{i_{5}}G_{i_{4}}G_{i_{3}}G_{i_{2}}B_{i_{1}}
$$

\vspace{3mm}
\noindent
are calculated at the point $({\bf y}_p,p).$

It is well known that
under the standard conditions \cite{KlPl2}, \cite{2006} the numerical 
scheme (\ref{4.470}) has 
strong order of convergence
2.5. The major emphasis below will be placed on the 
approximation of the iterated
Stratonovich stochastic integrals appearing in (\ref{4.470}). 
Therefore, among 
the standard conditions, we note the
following appro\-xi\-ma\-ti\-on condi\-ti\-on for these 
stochastic integrals \cite{KlPl2}, \cite{2006}

\vspace{-1mm}
\begin{equation}
\label{ors}
{\sf M}\biggl\{\biggl(I_{(l_{1}\ldots\hspace{0.2mm} l_{k})\tau_{p+1},\tau_p}
^{*(i_{1}\ldots i_{k})} 
- \hat I_{(l_{1}\ldots\hspace{0.2mm} l_{k})\tau_{p+1},
\tau_p}^{*(i_{1}\ldots i_{k})}
\biggr)^2\biggr\}\le C\Delta^{6},
\end{equation}

\vspace{3mm}
\noindent
where constant $C$ is independent of
$\Delta$.

Note that if we exclude from (\ref{4.470}) the terms starting from the
term $\Delta^3 LL{\bf a}/6$, then we have the explicit 
one-step strong numerical scheme with order 2.0 \cite{KlPl2}, 
\cite{2006}, \cite{2017-1}-\cite{2010-1}.

Using the numerical scheme (\ref{4.470}) or its modifications based 
on the classical Taylor--Stratonovich expansion \cite{KlPl1},
the implicit or multistep analogues of (\ref{4.470}) can be constructed
\cite{KlPl2}, \cite{2006}, \cite{2017-1}-\cite{2010-1}. The set of the
iterated Stratonovich 
stochastic integrals to be approximated for implementing 
these modifications is the same
as for the numerical scheme (\ref{4.470}) itself.
Interestingly, the truncated unified Taylor--Stratonovich expansion \cite{kk6}
(the 
foundation of the numerical
scheme (\ref{4.470})) contains only $12$ 
different types of the iterated Stratonovich
stochastic integrals 
(\ref{str11}), which cannot be
interconnected by linear relations \cite{2006}, \cite{2017-1}-\cite{2010-1}. 
The analogues 
classical Taylor--Stratonovich expansion \cite{KlPl2}, \cite{KlPl1} contains
$17$ different types of iterated Stratonovich 
stochastic integrals, part of which 
are interconnected by linear relations
and part of which have a higher multiplicity than the iterated 
Stratonovich stochastic integrals (\ref{str11}). This
fact well explains the use of the numerical scheme (\ref{4.470}).

One of the main problems arising in the implementation of the 
numerical scheme (\ref{4.470}) is the joint
numerical modeling of the iterated Stratonovich stochastic integrals 
figuring in (\ref{4.470}).

\vspace{5mm}

\section{Expansions 
of Iterated Ito 
Stochastic Integrals
(Method of Genegalized Multiple Fourier Series)}

\vspace{5mm}

An efficient numerical modeling method for iterated Ito 
stochastic integrals based on generalized multiple
Fourier series was considered in \cite{2006} (also see
\cite{2011-2}-\cite{200a},
\cite{301a}-\cite{new-art-1-xxy}).

This method rests on important
results presented below (Theorems 1, 2).

Suppose that every $\psi_l(\tau)$ $(l=1,\ldots,k)$ is a function from the space $L_2([t, T])$.
Define the following function on the hypercube $[t, T]^k$

\vspace{1mm}
\begin{equation}
\label{ppp}
K(t_1,\ldots,t_k)=
\begin{cases}
\psi_1(t_1)\ldots \psi_k(t_k)\ &\hbox{for}\ \ t_1<\ldots<t_k\\
~\\
~\\
0\ &\hbox{otherwise}
\end{cases},\ \ \ \ t_1,\ldots,t_k\in[t, T],\ \ \ \ k\ge 2,
\end{equation}

\vspace{5mm}
\noindent
and 
$K(t_1)\equiv\psi_1(t_1)$ for $t_1\in[t, T].$

Suppose that $\{\phi_j(x)\}_{j=0}^{\infty}$
is a complete orthonormal system of functions in the space
$L_2([t, T])$. 
The function $K(t_1,\ldots,t_k)$ belongs to the space $L_2([t, T]^k).$
At this situation it is well known that the generalized 
multiple Fourier series 
of $K(t_1,\ldots,t_k)\in L_2([t, T]^k)$ is converging 
to $K(t_1,\ldots,t_k)$ in the hypercube $[t, T]^k$ in 
the mean-square sense, i.e.

\vspace{2mm}
$$
\hbox{\vtop{\offinterlineskip\halign{
\hfil#\hfil\cr
{\rm lim}\cr
$\stackrel{}{{}_{p_1,\ldots,p_k\to \infty}}$\cr
}} }\Biggl\Vert
K(t_1,\ldots,t_k)-
\sum_{j_1=0}^{p_1}\ldots \sum_{j_k=0}^{p_k}
C_{j_k\ldots j_1}\prod_{l=1}^{k} \phi_{j_l}(t_l)
\Biggr\Vert_{L_2([t,T]^k)}=0,
$$

\vspace{5mm}
\noindent
where
\begin{equation}
\label{ppppa}
C_{j_k\ldots j_1}=\int\limits_{[t,T]^k}
K(t_1,\ldots,t_k)\prod_{l=1}^{k}\phi_{j_l}(t_l)dt_1\ldots dt_k
\end{equation}

\vspace{4mm}
\noindent
is the Fourier coefficient,

$$
\left\Vert f\right\Vert_{L_2([t,T]^k)}=\left(\int\limits_{[t,T]^k}
f^2(t_1,\ldots,t_k)dt_1\ldots dt_k\right)^{1/2}.
$$

\vspace{6mm}

Consider the partition $\{\tau_j\}_{j=0}^N$ of $[t,T]$ such that

\vspace{1mm}
\begin{equation}
\label{1111}
t=\tau_0<\ldots <\tau_N=T,\ \ \
\Delta_N=
\hbox{\vtop{\offinterlineskip\halign{
\hfil#\hfil\cr
{\rm max}\cr
$\stackrel{}{{}_{0\le j\le N-1}}$\cr
}} }\Delta\tau_j\to 0\ \ \hbox{if}\ \ N\to \infty,\ \ \
\Delta\tau_j=\tau_{j+1}-\tau_j.
\end{equation}

\vspace{4mm}

{\bf Theorem 1}\ \cite{2006} (2006), \cite{2011-2}-\cite{200a},
\cite{301a}-\cite{new-art-1-xxy}.\
{\it Suppose that
every $\psi_l(\tau)$ $(l=1,\ldots, k)$ is a continuous 
nonrandom func\-tion on 
$[t, T]$ and
$\{\phi_j(x)\}_{j=0}^{\infty}$ is a complete orthonormal system  
of continuous func\-ti\-ons in the space $L_2([t,T]).$ Then

\vspace{1mm}
$$
J[\psi^{(k)}]_{T,t}\  =\ 
\hbox{\vtop{\offinterlineskip\halign{
\hfil#\hfil\cr
{\rm l.i.m.}\cr
$\stackrel{}{{}_{p_1,\ldots,p_k\to \infty}}$\cr
}} }\sum_{j_1=0}^{p_1}\ldots\sum_{j_k=0}^{p_k}
C_{j_k\ldots j_1}\Biggl(
\prod_{l=1}^k\zeta_{j_l}^{(i_l)}\ -
\Biggr.
$$

\vspace{2mm}
\begin{equation}
\label{tyyy}
-\ \Biggl.
\hbox{\vtop{\offinterlineskip\halign{
\hfil#\hfil\cr
{\rm l.i.m.}\cr
$\stackrel{}{{}_{N\to \infty}}$\cr
}} }\sum_{(l_1,\ldots,l_k)\in {\rm G}_k}
\phi_{j_{1}}(\tau_{l_1})
\Delta{\bf w}_{\tau_{l_1}}^{(i_1)}\ldots
\phi_{j_{k}}(\tau_{l_k})
\Delta{\bf w}_{\tau_{l_k}}^{(i_k)}\Biggr),
\end{equation}

\vspace{5mm}
\noindent
where $J[\psi^{(k)}]_{T,t}$ is defined by {\rm (\ref{ito}),}

\vspace{1mm}
$$
{\rm G}_k={\rm H}_k\backslash{\rm L}_k,\ \ \
{\rm H}_k=\{(l_1,\ldots,l_k):\ l_1,\ldots,l_k=0,\ 1,\ldots,N-1\},
$$

\vspace{1mm}
$$
{\rm L}_k=\{(l_1,\ldots,l_k):\ l_1,\ldots,l_k=0,\ 1,\ldots,N-1;\
l_g\ne l_r\ (g\ne r);\ g, r=1,\ldots,k\},
$$

\vspace{5mm}
\noindent
${\rm l.i.m.}$ is a limit in the mean-square sense$,$
$i_1,\ldots,i_k=0,1,\ldots,m,$

\vspace{-1mm}
\begin{equation}
\label{rr23}
\zeta_{j}^{(i)}=
\int\limits_t^T \phi_{j}(s) d{\bf w}_s^{(i)}
\end{equation} 

\vspace{3mm}
\noindent
are independent standard Gaussian random variables
for various
$i$ or $j$ {\rm(}in the case when $i\ne 0${\rm),}
$C_{j_k\ldots j_1}$ is the Fourier coefficient {\rm(\ref{ppppa}),}
$\Delta{\bf w}_{\tau_{j}}^{(i)}=
{\bf w}_{\tau_{j+1}}^{(i)}-{\bf w}_{\tau_{j}}^{(i)}$
$(i=0, 1,\ldots,m),$
$\left\{\tau_{j}\right\}_{j=0}^{N}$ is a partition of
the interval $[t, T],$ which satisfies the condition {\rm (\ref{1111})}.
}

\vspace{2mm}

It was shown in \cite{2007-2}-\cite{2013}
that Theorem 1 is valid for convergence 
in the mean of degree $2n$ ($n\in \mathbb{N}$).
The convergence with probability 1 in Theorem 1  
is proved in \cite{2018a}-\cite{xxx333}, \cite{arxiv-1}, \cite{OK1000}
for the cases of Legendre polynomials and trigonometric functions.
Moreover, the complete orthonormal systems of Haar and 
Rademacher--Walsh functions in the space $L_2([t,T])$ 
can also be applied in Theorem 1
\cite{2006}-\cite{2013}.
The modification of Theorem 1 for 
complete orthonormal with weigth $r(x)\ge 0$ systems
of functions in the space $L_2([t,T])$ can be found in 
\cite{2018}, \cite{2018a}-\cite{xxx333}, \cite{arxiv-1}, \cite{arxiv-26b}.
Application of Theorem 1 and Theorem 2 (see below) for the mean-square
approximation of iterated stochastic integrals 
with respect to the 
infinite-dimensional $Q$-Wiener process can be found
in the monographs \cite{2018a}-\cite{xxx333} (Chapter 7) and in \cite{31a}, \cite{200a},
\cite{OK}, \cite{Kuzh-1}.

In order to evaluate the significance of Theorem 1 for practice we will
demonstrate its transfor\-med particular cases for 
$k=1,\ldots,5$ 
\cite{2006}-\cite{200a}, \cite{301a}-\cite{new-art-1-xxy}

\begin{equation}
\label{a1}
J[\psi^{(1)}]_{T,t}
=\hbox{\vtop{\offinterlineskip\halign{
\hfil#\hfil\cr
{\rm l.i.m.}\cr
$\stackrel{}{{}_{p_1\to \infty}}$\cr
}} }\sum_{j_1=0}^{p_1}
C_{j_1}\zeta_{j_1}^{(i_1)},
\end{equation}

\vspace{2mm}
\begin{equation}
\label{a2}
J[\psi^{(2)}]_{T,t}
=\hbox{\vtop{\offinterlineskip\halign{
\hfil#\hfil\cr
{\rm l.i.m.}\cr
$\stackrel{}{{}_{p_1,p_2\to \infty}}$\cr
}} }\sum_{j_1=0}^{p_1}\sum_{j_2=0}^{p_2}
C_{j_2j_1}\Biggl(\zeta_{j_1}^{(i_1)}\zeta_{j_2}^{(i_2)}
-{\bf 1}_{\{i_1=i_2\ne 0\}}
{\bf 1}_{\{j_1=j_2\}}\Biggr),
\end{equation}

\vspace{4mm}

$$
J[\psi^{(3)}]_{T,t}=
\hbox{\vtop{\offinterlineskip\halign{
\hfil#\hfil\cr
{\rm l.i.m.}\cr
$\stackrel{}{{}_{p_1,\ldots,p_3\to \infty}}$\cr
}} }\sum_{j_1=0}^{p_1}\sum_{j_2=0}^{p_2}\sum_{j_3=0}^{p_3}
C_{j_3j_2j_1}\Biggl(
\zeta_{j_1}^{(i_1)}\zeta_{j_2}^{(i_2)}\zeta_{j_3}^{(i_3)}
-\Biggr.
$$

\begin{equation}
\label{a3}
-\Biggl.
{\bf 1}_{\{i_1=i_2\ne 0\}}
{\bf 1}_{\{j_1=j_2\}}
\zeta_{j_3}^{(i_3)}
-{\bf 1}_{\{i_2=i_3\ne 0\}}
{\bf 1}_{\{j_2=j_3\}}
\zeta_{j_1}^{(i_1)}-
{\bf 1}_{\{i_1=i_3\ne 0\}}
{\bf 1}_{\{j_1=j_3\}}
\zeta_{j_2}^{(i_2)}\Biggr),
\end{equation}

\vspace{6mm}

$$
J[\psi^{(4)}]_{T,t}
=
\hbox{\vtop{\offinterlineskip\halign{
\hfil#\hfil\cr
{\rm l.i.m.}\cr
$\stackrel{}{{}_{p_1,\ldots,p_4\to \infty}}$\cr
}} }\sum_{j_1=0}^{p_1}\ldots\sum_{j_4=0}^{p_4}
C_{j_4\ldots j_1}\Biggl(
\prod_{l=1}^4\zeta_{j_l}^{(i_l)}
\Biggr.
-
$$
$$
-
{\bf 1}_{\{i_1=i_2\ne 0\}}
{\bf 1}_{\{j_1=j_2\}}
\zeta_{j_3}^{(i_3)}
\zeta_{j_4}^{(i_4)}
-
{\bf 1}_{\{i_1=i_3\ne 0\}}
{\bf 1}_{\{j_1=j_3\}}
\zeta_{j_2}^{(i_2)}
\zeta_{j_4}^{(i_4)}-
$$
$$
-
{\bf 1}_{\{i_1=i_4\ne 0\}}
{\bf 1}_{\{j_1=j_4\}}
\zeta_{j_2}^{(i_2)}
\zeta_{j_3}^{(i_3)}
-
{\bf 1}_{\{i_2=i_3\ne 0\}}
{\bf 1}_{\{j_2=j_3\}}
\zeta_{j_1}^{(i_1)}
\zeta_{j_4}^{(i_4)}-
$$
$$
-
{\bf 1}_{\{i_2=i_4\ne 0\}}
{\bf 1}_{\{j_2=j_4\}}
\zeta_{j_1}^{(i_1)}
\zeta_{j_3}^{(i_3)}
-
{\bf 1}_{\{i_3=i_4\ne 0\}}
{\bf 1}_{\{j_3=j_4\}}
\zeta_{j_1}^{(i_1)}
\zeta_{j_2}^{(i_2)}+
$$
$$
+
{\bf 1}_{\{i_1=i_2\ne 0\}}
{\bf 1}_{\{j_1=j_2\}}
{\bf 1}_{\{i_3=i_4\ne 0\}}
{\bf 1}_{\{j_3=j_4\}}
+
$$
$$
+
{\bf 1}_{\{i_1=i_3\ne 0\}}
{\bf 1}_{\{j_1=j_3\}}
{\bf 1}_{\{i_2=i_4\ne 0\}}
{\bf 1}_{\{j_2=j_4\}}+
$$
\begin{equation}
\label{a4}
+\Biggl.
{\bf 1}_{\{i_1=i_4\ne 0\}}
{\bf 1}_{\{j_1=j_4\}}
{\bf 1}_{\{i_2=i_3\ne 0\}}
{\bf 1}_{\{j_2=j_3\}}\Biggr),
\end{equation}

\vspace{7mm}

$$
J[\psi^{(5)}]_{T,t}
=\hbox{\vtop{\offinterlineskip\halign{
\hfil#\hfil\cr
{\rm l.i.m.}\cr
$\stackrel{}{{}_{p_1,\ldots,p_5\to \infty}}$\cr
}} }\sum_{j_1=0}^{p_1}\ldots\sum_{j_5=0}^{p_5}
C_{j_5\ldots j_1}\Biggl(
\prod_{l=1}^5\zeta_{j_l}^{(i_l)}
-\Biggr.
$$
$$
-
{\bf 1}_{\{i_1=i_2\ne 0\}}
{\bf 1}_{\{j_1=j_2\}}
\zeta_{j_3}^{(i_3)}
\zeta_{j_4}^{(i_4)}
\zeta_{j_5}^{(i_5)}-
{\bf 1}_{\{i_1=i_3\ne 0\}}
{\bf 1}_{\{j_1=j_3\}}
\zeta_{j_2}^{(i_2)}
\zeta_{j_4}^{(i_4)}
\zeta_{j_5}^{(i_5)}-
$$
$$
-
{\bf 1}_{\{i_1=i_4\ne 0\}}
{\bf 1}_{\{j_1=j_4\}}
\zeta_{j_2}^{(i_2)}
\zeta_{j_3}^{(i_3)}
\zeta_{j_5}^{(i_5)}-
{\bf 1}_{\{i_1=i_5\ne 0\}}
{\bf 1}_{\{j_1=j_5\}}
\zeta_{j_2}^{(i_2)}
\zeta_{j_3}^{(i_3)}
\zeta_{j_4}^{(i_4)}-
$$
$$
-
{\bf 1}_{\{i_2=i_3\ne 0\}}
{\bf 1}_{\{j_2=j_3\}}
\zeta_{j_1}^{(i_1)}
\zeta_{j_4}^{(i_4)}
\zeta_{j_5}^{(i_5)}-
{\bf 1}_{\{i_2=i_4\ne 0\}}
{\bf 1}_{\{j_2=j_4\}}
\zeta_{j_1}^{(i_1)}
\zeta_{j_3}^{(i_3)}
\zeta_{j_5}^{(i_5)}-
$$
$$
-
{\bf 1}_{\{i_2=i_5\ne 0\}}
{\bf 1}_{\{j_2=j_5\}}
\zeta_{j_1}^{(i_1)}
\zeta_{j_3}^{(i_3)}
\zeta_{j_4}^{(i_4)}
-{\bf 1}_{\{i_3=i_4\ne 0\}}
{\bf 1}_{\{j_3=j_4\}}
\zeta_{j_1}^{(i_1)}
\zeta_{j_2}^{(i_2)}
\zeta_{j_5}^{(i_5)}-
$$
$$
-
{\bf 1}_{\{i_3=i_5\ne 0\}}
{\bf 1}_{\{j_3=j_5\}}
\zeta_{j_1}^{(i_1)}
\zeta_{j_2}^{(i_2)}
\zeta_{j_4}^{(i_4)}
-{\bf 1}_{\{i_4=i_5\ne 0\}}
{\bf 1}_{\{j_4=j_5\}}
\zeta_{j_1}^{(i_1)}
\zeta_{j_2}^{(i_2)}
\zeta_{j_3}^{(i_3)}+
$$
$$
+
{\bf 1}_{\{i_1=i_2\ne 0\}}
{\bf 1}_{\{j_1=j_2\}}
{\bf 1}_{\{i_3=i_4\ne 0\}}
{\bf 1}_{\{j_3=j_4\}}\zeta_{j_5}^{(i_5)}+
{\bf 1}_{\{i_1=i_2\ne 0\}}
{\bf 1}_{\{j_1=j_2\}}
{\bf 1}_{\{i_3=i_5\ne 0\}}
{\bf 1}_{\{j_3=j_5\}}\zeta_{j_4}^{(i_4)}+
$$
$$
+
{\bf 1}_{\{i_1=i_2\ne 0\}}
{\bf 1}_{\{j_1=j_2\}}
{\bf 1}_{\{i_4=i_5\ne 0\}}
{\bf 1}_{\{j_4=j_5\}}\zeta_{j_3}^{(i_3)}+
{\bf 1}_{\{i_1=i_3\ne 0\}}
{\bf 1}_{\{j_1=j_3\}}
{\bf 1}_{\{i_2=i_4\ne 0\}}
{\bf 1}_{\{j_2=j_4\}}\zeta_{j_5}^{(i_5)}+
$$
$$
+
{\bf 1}_{\{i_1=i_3\ne 0\}}
{\bf 1}_{\{j_1=j_3\}}
{\bf 1}_{\{i_2=i_5\ne 0\}}
{\bf 1}_{\{j_2=j_5\}}\zeta_{j_4}^{(i_4)}+
{\bf 1}_{\{i_1=i_3\ne 0\}}
{\bf 1}_{\{j_1=j_3\}}
{\bf 1}_{\{i_4=i_5\ne 0\}}
{\bf 1}_{\{j_4=j_5\}}\zeta_{j_2}^{(i_2)}+
$$
$$
+
{\bf 1}_{\{i_1=i_4\ne 0\}}
{\bf 1}_{\{j_1=j_4\}}
{\bf 1}_{\{i_2=i_3\ne 0\}}
{\bf 1}_{\{j_2=j_3\}}\zeta_{j_5}^{(i_5)}+
{\bf 1}_{\{i_1=i_4\ne 0\}}
{\bf 1}_{\{j_1=j_4\}}
{\bf 1}_{\{i_2=i_5\ne 0\}}
{\bf 1}_{\{j_2=j_5\}}\zeta_{j_3}^{(i_3)}+
$$
$$
+
{\bf 1}_{\{i_1=i_4\ne 0\}}
{\bf 1}_{\{j_1=j_4\}}
{\bf 1}_{\{i_3=i_5\ne 0\}}
{\bf 1}_{\{j_3=j_5\}}\zeta_{j_2}^{(i_2)}+
{\bf 1}_{\{i_1=i_5\ne 0\}}
{\bf 1}_{\{j_1=j_5\}}
{\bf 1}_{\{i_2=i_3\ne 0\}}
{\bf 1}_{\{j_2=j_3\}}\zeta_{j_4}^{(i_4)}+
$$
$$
+
{\bf 1}_{\{i_1=i_5\ne 0\}}
{\bf 1}_{\{j_1=j_5\}}
{\bf 1}_{\{i_2=i_4\ne 0\}}
{\bf 1}_{\{j_2=j_4\}}\zeta_{j_3}^{(i_3)}+
{\bf 1}_{\{i_1=i_5\ne 0\}}
{\bf 1}_{\{j_1=j_5\}}
{\bf 1}_{\{i_3=i_4\ne 0\}}
{\bf 1}_{\{j_3=j_4\}}\zeta_{j_2}^{(i_2)}+
$$
$$
+
{\bf 1}_{\{i_2=i_3\ne 0\}}
{\bf 1}_{\{j_2=j_3\}}
{\bf 1}_{\{i_4=i_5\ne 0\}}
{\bf 1}_{\{j_4=j_5\}}\zeta_{j_1}^{(i_1)}+
{\bf 1}_{\{i_2=i_4\ne 0\}}
{\bf 1}_{\{j_2=j_4\}}
{\bf 1}_{\{i_3=i_5\ne 0\}}
{\bf 1}_{\{j_3=j_5\}}\zeta_{j_1}^{(i_1)}+
$$
\begin{equation}
\label{a5}
+\Biggl.
{\bf 1}_{\{i_2=i_5\ne 0\}}
{\bf 1}_{\{j_2=j_5\}}
{\bf 1}_{\{i_3=i_4\ne 0\}}
{\bf 1}_{\{j_3=j_4\}}\zeta_{j_1}^{(i_1)}\Biggr),
\end{equation}

\vspace{6mm}
\noindent
where ${\bf 1}_A$ is the indicator of the set $A$.

Note that we will consider the case $i_1,\ldots,i_5=1,\ldots,m$.
This case corresponds to the numerical scheme (\ref{4.470}).

For further consideration, let us 
consider the generalization of formulas (\ref{a1})--(\ref{a5})                 
for the case of an arbitrary multiplicity $k$ $(k\in\mathbb{N})$ of 
the iterated Ito stochastic integral $J[\psi^{(k)}]_{T,t}$ defined by (\ref{ito}).
In order to do this, let us
introduce some notations. 
Consider the unordered
set $\{1, 2, \ldots, k\}$ 
and separate it into two parts:
the first part consists of $r$ unordered 
pairs (sequence order of these pairs is also unimportant) and the 
second one consists of the 
remaining $k-2r$ numbers.
So, we have

\begin{equation}
\label{leto5007}
(\{
\underbrace{\{g_1, g_2\}, \ldots, 
\{g_{2r-1}, g_{2r}\}}_{\small{\hbox{part 1}}}
\},
\{\underbrace{q_1, \ldots, q_{k-2r}}_{\small{\hbox{part 2}}}
\}),
\end{equation}

\vspace{4mm}
\noindent
where 

\vspace{-2mm}
$$
\{g_1, g_2, \ldots, 
g_{2r-1}, g_{2r}, q_1, \ldots, q_{k-2r}\}=\{1, 2, \ldots, k\},
$$

\vspace{4mm}
\noindent
braces   
mean an unordered 
set, and pa\-ren\-the\-ses mean an ordered set.

We will say that (\ref{leto5007}) is a partition 
and consider the sum with respect to all possible
partitions

\begin{equation}
\label{leto5008}
\sum_{\stackrel{(\{\{g_1, g_2\}, \ldots, 
\{g_{2r-1}, g_{2r}\}\}, \{q_1, \ldots, q_{k-2r}\})}
{{}_{\{g_1, g_2, \ldots, 
g_{2r-1}, g_{2r}, q_1, \ldots, q_{k-2r}\}=\{1, 2, \ldots, k\}}}}
a_{g_1 g_2, \ldots, 
g_{2r-1} g_{2r}, q_1 \ldots q_{k-2r}}.
\end{equation}

\vspace{4mm}

Below there are several examples of sums in the form (\ref{leto5008})

\vspace{2mm}
$$
\sum_{\stackrel{(\{g_1, g_2\})}{{}_{\{g_1, g_2\}=\{1, 2\}}}}
a_{g_1 g_2}=a_{12},
$$

\vspace{3mm}
$$
\sum_{\stackrel{(\{\{g_1, g_2\}, \{g_3, g_4\}\})}
{{}_{\{g_1, g_2, g_3, g_4\}=\{1, 2, 3, 4\}}}}
a_{g_1 g_2 g_3 g_4}=a_{1234} + a_{1324} + a_{2314},
$$

\vspace{3mm}
$$
\sum_{\stackrel{(\{g_1, g_2\}, \{q_1, q_{2}\})}
{{}_{\{g_1, g_2, q_1, q_{2}\}=\{1, 2, 3, 4\}}}}
a_{g_1 g_2, q_1 q_{2}}=
$$

$$
=a_{12,34}+a_{13,24}+a_{14,23}
+a_{23,14}+a_{24,13}+a_{34,12},
$$

\vspace{3mm}
$$
\sum_{\stackrel{(\{g_1, g_2\}, \{q_1, q_{2}, q_3\})}
{{}_{\{g_1, g_2, q_1, q_{2}, q_3\}=\{1, 2, 3, 4, 5\}}}}
a_{g_1 g_2, q_1 q_{2}q_3}
=
$$

$$
=a_{12,345}+a_{13,245}+a_{14,235}
+a_{15,234}+a_{23,145}+a_{24,135}+
$$
$$
+a_{25,134}+a_{34,125}+a_{35,124}+a_{45,123},
$$

\vspace{4mm}
$$
\sum_{\stackrel{(\{\{g_1, g_2\}, \{g_3, g_{4}\}\}, \{q_1\})}
{{}_{\{g_1, g_2, g_3, g_{4}, q_1\}=\{1, 2, 3, 4, 5\}}}}
a_{g_1 g_2, g_3 g_{4},q_1}
=
$$

$$
=
a_{12,34,5}+a_{13,24,5}+a_{14,23,5}+
a_{12,35,4}+a_{13,25,4}+a_{15,23,4}+
$$
$$
+a_{12,54,3}+a_{15,24,3}+a_{14,25,3}+a_{15,34,2}+a_{13,54,2}+a_{14,53,2}+
$$
$$
+
a_{52,34,1}+a_{53,24,1}+a_{54,23,1}.
$$

\vspace{5mm}

Now we can write (\ref{tyyy}) as

\vspace{1mm}

$$
J[\psi^{(k)}]_{T,t}=
\hbox{\vtop{\offinterlineskip\halign{
\hfil#\hfil\cr
{\rm l.i.m.}\cr
$\stackrel{}{{}_{p_1,\ldots,p_k\to \infty}}$\cr
}} }
\sum\limits_{j_1=0}^{p_1}\ldots
\sum\limits_{j_k=0}^{p_k}
C_{j_k\ldots j_1}\Biggl(
\prod_{l=1}^k\zeta_{j_l}^{(i_l)}+\sum\limits_{r=1}^{[k/2]}
(-1)^r \times
\Biggr.
$$

\vspace{3mm}
\begin{equation}
\label{leto6000hh}
\times
\sum_{\stackrel{(\{\{g_1, g_2\}, \ldots, 
\{g_{2r-1}, g_{2r}\}\}, \{q_1, \ldots, q_{k-2r}\})}
{{}_{\{g_1, g_2, \ldots, 
g_{2r-1}, g_{2r}, q_1, \ldots, q_{k-2r}\}=\{1, 2, \ldots, k\}}}}
\prod\limits_{s=1}^r
{\bf 1}_{\{i_{g_{{}_{2s-1}}}=~i_{g_{{}_{2s}}}\ne 0\}}
\Biggl.{\bf 1}_{\{j_{g_{{}_{2s-1}}}=~j_{g_{{}_{2s}}}\}}
\prod_{l=1}^{k-2r}\zeta_{j_{q_l}}^{(i_{q_l})}\Biggr),
\end{equation}

\vspace{5mm}
\noindent
where $[x]$ is an integer part of a real number $x;$
another notations are the same as in Theorem {\bf 1}.

\vspace{2mm}

In particular, from (\ref{leto6000hh}) for $k=5$ we obtain

\vspace{3mm}

$$
J[\psi^{(5)}]_{T,t}=
\hbox{\vtop{\offinterlineskip\halign{
\hfil#\hfil\cr
{\rm l.i.m.}\cr
$\stackrel{}{{}_{p_1,\ldots,p_5\to \infty}}$\cr
}} }\sum_{j_1=0}^{p_1}\ldots\sum_{j_5=0}^{p_5}
C_{j_5\ldots j_1}\Biggl(
\prod_{l=1}^5\zeta_{j_l}^{(i_l)}-\Biggr.
$$

\vspace{2mm}
$$
-
\sum\limits_{\stackrel{(\{g_1, g_2\}, \{q_1, q_{2}, q_3\})}
{{}_{\{g_1, g_2, q_{1}, q_{2}, q_3\}=\{1, 2, 3, 4, 5\}}}}
{\bf 1}_{\{i_{g_{{}_{1}}}=~i_{g_{{}_{2}}}\ne 0\}}
{\bf 1}_{\{j_{g_{{}_{1}}}=~j_{g_{{}_{2}}}\}}
\prod_{l=1}^{3}\zeta_{j_{q_l}}^{(i_{q_l})}+
$$

\vspace{2mm}
$$
+
\sum_{\stackrel{(\{\{g_1, g_2\}, 
\{g_{3}, g_{4}\}\}, \{q_1\})}
{{}_{\{g_1, g_2, g_{3}, g_{4}, q_1\}=\{1, 2, 3, 4, 5\}}}}
{\bf 1}_{\{i_{g_{{}_{1}}}=~i_{g_{{}_{2}}}\ne 0\}}
{\bf 1}_{\{j_{g_{{}_{1}}}=~j_{g_{{}_{2}}}\}}
\Biggl.{\bf 1}_{\{i_{g_{{}_{3}}}=~i_{g_{{}_{4}}}\ne 0\}}
{\bf 1}_{\{j_{g_{{}_{3}}}=~j_{g_{{}_{4}}}\}}
\zeta_{j_{q_1}}^{(i_{q_1})}\Biggr).
$$

\vspace{7mm}
\noindent
The last equality obviously agrees with
(\ref{a5}).

Let us consider a generalization of Theorem 1 for the case
of an arbitrary complete orthonormal systems  
of functions in the space $L_2([t,T])$ 
and $\psi_1(\tau),\ldots,\psi_k(\tau)\in L_2([t, T]).$

\vspace{2mm}

{\bf Theorem~2}\ \cite{2018a} (Sect.~1.11), \cite{arxiv-1} (Sect.~15).
{\it Suppose that
$\psi_1(\tau),\ldots,\psi_k(\tau)\in L_2([t, T])$ and
$\{\phi_j(x)\}_{j=0}^{\infty}$ is an arbitrary complete orthonormal system  
of functions in the space $L_2([t,T]).$
Then the following expansion

\vspace{1mm}
$$
J[\psi^{(k)}]_{T,t}=
\hbox{\vtop{\offinterlineskip\halign{
\hfil#\hfil\cr
{\rm l.i.m.}\cr
$\stackrel{}{{}_{p_1,\ldots,p_k\to \infty}}$\cr
}} }
\sum\limits_{j_1=0}^{p_1}\ldots
\sum\limits_{j_k=0}^{p_k}
C_{j_k\ldots j_1}\Biggl(
\prod_{l=1}^k\zeta_{j_l}^{(i_l)}+\sum\limits_{r=1}^{[k/2]}
(-1)^r \times
\Biggr.
$$

\vspace{2mm}
\begin{equation}
\label{leto6000}
\times
\sum_{\stackrel{(\{\{g_1, g_2\}, \ldots, 
\{g_{2r-1}, g_{2r}\}\}, \{q_1, \ldots, q_{k-2r}\})}
{{}_{\{g_1, g_2, \ldots, 
g_{2r-1}, g_{2r}, q_1, \ldots, q_{k-2r}\}=\{1, 2, \ldots, k\}}}}
\prod\limits_{s=1}^r
{\bf 1}_{\{i_{g_{{}_{2s-1}}}=~i_{g_{{}_{2s}}}\ne 0\}}
\Biggl.{\bf 1}_{\{j_{g_{{}_{2s-1}}}=~j_{g_{{}_{2s}}}\}}
\prod_{l=1}^{k-2r}\zeta_{j_{q_l}}^{(i_{q_l})}\Biggr)
\end{equation}

\vspace{6mm}
\noindent
con\-verg\-ing in the mean-square sense is valid,
where $[x]$ is an integer part of a real number $x;$
another notations are the same as in Theorem~{\rm 1}.}

\vspace{2mm}

It should be noted that an analogue of Theorem 2 was considered 
in \cite{Rybakov1000}. 
Note that we use another notations 
\cite{2018a} (Sect.~1.11), \cite{arxiv-1} (Sect.~15)
in comparison with \cite{Rybakov1000}.
Moreover, the proof of an analogue of Theorem 2
from \cite{Rybakov1000} is somewhat different from the proof given in 
\cite{2018a} (Sect.~1.11), \cite{arxiv-1} (Sect.~15).

\vspace{5mm}

\section{Calculation of the Mean-Square Approximation Error
in the Method of Generalized Multiple Fourier Seires}

\vspace{5mm}

Note that for the integrals $J[\psi^{(k)}]_{T,t}$ defined by 
(\ref{ito})
the mean-square approximation error can be exactly
calculated and efficiently estimated.

Let $J[\psi^{(k)}]_{T,t}^{q}$ be the
expression on the right-hand side of (\ref{leto6000}) before passing to the limit
$\hbox{\vtop{\offinterlineskip\halign{
\hfil#\hfil\cr
{\rm l.i.m.}\cr
$\stackrel{}{{}_{p_1,\ldots,p_k\to \infty}}$\cr
}} }$ for the case
$p_1=\ldots=p_k=q,$ i.e.

\vspace{1mm}

$$
J[\psi^{(k)}]_{T,t}^{q}=
\sum\limits_{j_1,\ldots,j_k=0}^{q}
C_{j_k\ldots j_1}\Biggl(
\prod_{l=1}^k\zeta_{j_l}^{(i_l)}+\sum\limits_{r=1}^{[k/2]}
(-1)^r \times
\Biggr.
$$

\vspace{2mm}
\begin{equation}
\label{r1}
\times
\sum_{\stackrel{(\{\{g_1, g_2\}, \ldots, 
\{g_{2r-1}, g_{2r}\}\}, \{q_1, \ldots, q_{k-2r}\})}
{{}_{\{g_1, g_2, \ldots, 
g_{2r-1}, g_{2r}, q_1, \ldots, q_{k-2r}\}=\{1, 2, \ldots, k\}}}}
\prod\limits_{s=1}^r
{\bf 1}_{\{i_{g_{{}_{2s-1}}}=~i_{g_{{}_{2s}}}\ne 0\}}
\Biggl.{\bf 1}_{\{j_{g_{{}_{2s-1}}}=~j_{g_{{}_{2s}}}\}}
\prod_{l=1}^{k-2r}\zeta_{j_{q_l}}^{(i_{q_l})}\Biggr).
\end{equation}

\vspace{5mm}

Let us denote

$$
{\sf M}\left\{\left(J[\psi^{(k)}]_{T,t}-
J[\psi^{(k)}]_{T,t}^{q}\right)^2\right\}\stackrel{{\rm def}}
{=}E_k^{q},
$$

\vspace{3mm}
$$
\int\limits_{[t,T]^k}
K^2(t_1,\ldots,t_k)dt_1\ldots dt_k
\stackrel{{\rm def}}{=}I_k.
$$

\vspace{6mm}

In \cite{2017-1}-\cite{xxx333}, \cite{arxiv-1}, 
\cite{arxiv-2} it was shown that

\begin{equation}
\label{qq4}
E_k^{q}\le k!\Biggl(I_k-\sum_{j_1,\ldots,j_k=0}^{q}C^2_{j_k\ldots j_1}\Biggr)
\end{equation}

\vspace{4mm}
\noindent
for the following two cases:

\vspace{2mm}

1.\ $i_1,\ldots,i_k=1,\ldots,m$ and $T-t\in (0, +\infty)$,

2.\ $i_1,\ldots,i_k=0, 1,\ldots,m$ and  $T-t\in (0, 1)$.

\vspace{3mm}

The value $E_k^{q}$
can be calculated exactly.

\vspace{2mm}

{\bf Theorem 3} \cite{2018a} (Sect.~1.12), \cite{arxiv-2} (Sect.~6).
{\it Suppose that $\{\phi_j(x)\}_{j=0}^{\infty}$ 
is an arbitrary complete orthonormal system  
of functions in the space $L_2([t,T])$ and
$\psi_1(\tau),\ldots,\psi_k(\tau)\in L_2([t, T]).$  
Then

\begin{equation}
\label{tttr11}
E_k^q=I_k- \sum_{j_1,\ldots, j_k=0}^{q}
C_{j_k\ldots j_1}
{\sf M}\left\{J[\psi^{(k)}]_{T,t}
\sum\limits_{(j_1,\ldots,j_k)}
\int\limits_t^T \phi_{j_k}(t_k)
\ldots
\int\limits_t^{t_{2}}\phi_{j_{1}}(t_{1})
d{\bf f}_{t_1}^{(i_1)}\ldots
d{\bf f}_{t_k}^{(i_k)}\right\},
\end{equation}

\vspace{5mm}
\noindent
where
$i_1,\ldots,i_k = 1,\ldots,m;$\ 
the expression 

\vspace{-1mm}
$$
\sum\limits_{(j_1,\ldots,j_k)}
$$ 

\vspace{3mm}
\noindent
means the sum with respect to all
possible permutations
$(j_1,\ldots,j_k)$. At the same time if 
$j_r$ swapped with $j_q$ in the permutation $(j_1,\ldots,j_k)$,
then $i_r$ swapped with $i_q$ in the permutation 
$(i_1,\ldots,i_k);$
another notations are the same as in Theorems {\rm 1, 2.}
}

\vspace{2mm}

Note that 

\vspace{-2mm}
$$
{\sf M}\left\{J[\psi^{(k)}]_{T,t}
\int\limits_t^T \phi_{j_k}(t_k)
\ldots
\int\limits_t^{t_{2}}\phi_{j_{1}}(t_{1})
d{\bf f}_{t_1}^{(i_1)}\ldots
d{\bf f}_{t_k}^{(i_k)}\right\}=C_{j_k\ldots j_1}.
$$

\vspace{4mm}

Therefore, for the case of pairwise 
different numbers $i_1,\ldots,i_k$ 
as well as for the case $i_1=\ldots=i_k$
from Theorem 3 it follows that
\cite{2018}, \cite{2018a}-\cite{xxx333}, \cite{17a}, 
\cite{arxiv-2}

\vspace{-1mm}
\begin{equation}
\label{qq1}
E_k^q= I_k- \sum_{j_1,\ldots,j_k=0}^{q}
C_{j_k\ldots j_1}^2,
\end{equation}

$$
E_k^q= I_k - \sum_{j_1,\ldots,j_k=0}^{q}
C_{j_k\ldots j_1}\Biggl(\sum\limits_{(j_1,\ldots,j_k)}
C_{j_k\ldots j_1}\Biggr),
$$

\vspace{4mm}
\noindent
where 

\vspace{-2mm}
$$
\sum\limits_{(j_1,\ldots,j_k)}
$$ 

\vspace{3mm}
\noindent
is a sum with respect to all 
possible permutations
$(j_1,\ldots,j_k)$.

Consider some examples \cite{2018}, \cite{2018a}-\cite{xxx333}, \cite{17a}, 
\cite{arxiv-2} of application of Theorem 3
$(i_1,i_2,i_3=1,\ldots,m)$

\vspace{1mm}
\begin{equation}
\label{qq2}
E_2^q     
=I_2
-\sum_{j_1,j_2=0}^q
C_{j_2j_1}^2-
\sum_{j_1,j_2=0}^q
C_{j_2j_1}C_{j_1j_2}\ \ \ (i_1=i_2),
\end{equation}

\vspace{2mm}
\begin{equation}
\label{qq3}
E_3^q=I_3
-\sum_{j_3,j_2,j_1=0}^q C_{j_3j_2j_1}^2-
\sum_{j_3,j_2,j_1=0}^q C_{j_3j_1j_2}C_{j_3j_2j_1}\ \ \ (i_1=i_2\ne i_3),
\end{equation}

\vspace{2mm}
\begin{equation}
\label{882}
E_3^q=I_3-
\sum_{j_3,j_2,j_1=0}^q C_{j_3j_2j_1}^2-
\sum_{j_3,j_2,j_1=0}^q C_{j_2j_3j_1}C_{j_3j_2j_1}\ \ \ (i_1\ne i_2=i_3),
\end{equation}

\vspace{2mm}
\begin{equation}
\label{883}
E_3^q=I_3
-\sum_{j_3,j_2,j_1=0}^q C_{j_3j_2j_1}^2-
\sum_{j_3,j_2,j_1=0}^q C_{j_3j_2j_1}C_{j_1j_2j_3}\ \ \ (i_1=i_3\ne i_2).
\end{equation}

\vspace{6mm}

The values $E_4^q$ and $E_5^q$ were calculated exaclty for all possible 
combinations of $i_1,\ldots,i_5=1,\ldots,m$ in 
\cite{2018}, \cite{2018a}-\cite{xxx333},  
\cite{arxiv-2}.

\vspace{5mm}

\section{Expansions 
of Iterated Stratonovich
Stochastic Integrals Based on Multiple Fourier--Legendre
Series and Multiple Trigonometric Fourier Seires}

\vspace{5mm}

In contrast to the iterated Ito stochastic integrals (\ref{ito}), 
the iterated Stratonovich stochastic integrals (\ref{str})
have simpler expansions (see Theorems 4--10 below) 
than (\ref{tyyy}) but the calculation (or estimation) 
of mean-square approximation
errors for the latter is a more difficult problem than 
for the former. We will study this issue in
details below.

As we mentioned above, Theorems 1, 2 can be adapted for the iterated
Stratonovich stochastic integrals (\ref{str}) at least
for multiplicities 1 to 6. 
Expansions of these iterated Stratonovich 
stochastic integrals turned out
much simpler, than the appropriate expansions
of the iterated Ito stochastic integrals (\ref{ito}) from Theorems 1, 2.
Let us formulate some old results on expansions of the iterated
Stratonovich stochastic integrals (\ref{str}) of
multiplicities 2 to 4.

\vspace{2mm}

{\bf Theorem 4}\ \cite{2011-2}-\cite{xxx333}, 
\cite{2010-2}-\cite{2013}, \cite{30a}, \cite{300a},
\cite{400a}, \cite{271a},
\cite{arxiv-5}, \cite{arxiv-8}, \cite{arxiv-23}.\
{\it Assume that the following conditions are fulfulled{\rm :}

{\rm 1}. The function $\psi_2(\tau)$ is continuously 
differentiable at the interval $[t, T]$ and the
function $\psi_1(\tau)$ is twice continuously 
differentiable at the interval $[t, T]$.

{\rm 2}. $\{\phi_j(x)\}_{j=0}^{\infty}$ is a complete orthonormal system 
of Legendre polynomials or tri\-go\-no\-met\-ric functions
in the space $L_2([t, T]).$

Then, the iterated Stratonovich stochastic integral of second multiplicity

$$
{\int\limits_t^{*}}^T\psi_2(t_2)
{\int\limits_t^{*}}^{t_2}\psi_1(t_1)d{\bf f}_{t_1}^{(i_1)}
d{\bf f}_{t_2}^{(i_2)}\ \ \ (i_1, i_2=1,\ldots,m)
$$

\vspace{4mm}
\noindent
is expanded into the 
following series

$$
{\int\limits_t^{*}}^T\psi_2(t_2)
{\int\limits_t^{*}}^{t_2}\psi_1(t_1)d{\bf f}_{t_1}^{(i_1)}=
\hbox{\vtop{\offinterlineskip\halign{
\hfil#\hfil\cr
{\rm l.i.m.}\cr
$\stackrel{}{{}_{p_1,p_2\to \infty}}$\cr
}} }\sum_{j_1=0}^{p_1}\sum_{j_2=0}^{p_2}
C_{j_2j_1}\zeta_{j_1}^{(i_1)}\zeta_{j_2}^{(i_2)}
$$

\vspace{4mm}
\noindent
converging in the mean-square sense, where 

$$
C_{j_2 j_1}=\int\limits_t^T
\psi_2(t_2)\phi_{j_2}(t_2)
\int\limits_t^{t_2}\psi_1(t_1)\phi_{j_1}(t_1)dt_1dt_2;
$$

\vspace{4mm}
\noindent
another notations are the same as in Theorems {\rm 1, 2}.}

\vspace{2mm}

{\bf Theorem 5}\ \cite{2011-2}-\cite{xxx333}, 
\cite{2010-2}-\cite{2013}, \cite{271a},
\cite{arxiv-5}, \cite{arxiv-7}.\
{\it Assume that
$\{\phi_j(x)\}_{j=0}^{\infty}$ is a complete orthonormal
system of Legendre polynomials or trigonomertic functions
in the space $L_2([t, T])$. Moreover,
the function $\psi_2(\tau)$ is continuously
differentiable at the interval $[t, T]$ and
the functions $\psi_1(\tau),$ $\psi_3(\tau)$ are twice continuously
differentiable at the interval $[t, T]$.

Then, for the iterated Stratonovich stochastic integral of 
third multiplicity

$$
{\int\limits_t^{*}}^T\psi_3(t_3)
{\int\limits_t^{*}}^{t_3}\psi_2(t_2)
{\int\limits_t^{*}}^{t_2}\psi_1(t_1)
d{\bf f}_{t_1}^{(i_1)}
d{\bf f}_{t_2}^{(i_2)}d{\bf f}_{t_3}^{(i_3)}\ \ \ (i_1, i_2, i_3=1,\ldots,m)
$$

\vspace{4mm}
\noindent
the following 
expansion 

\begin{equation}
\label{feto19000a}
{\int\limits_t^{*}}^T\psi_3(t_3)
{\int\limits_t^{*}}^{t_3}\psi_2(t_2)
{\int\limits_t^{*}}^{t_2}\psi_1(t_1)
d{\bf f}_{t_1}^{(i_1)}
d{\bf f}_{t_2}^{(i_2)}d{\bf f}_{t_3}^{(i_3)}\ 
=
\hbox{\vtop{\offinterlineskip\halign{
\hfil#\hfil\cr
{\rm l.i.m.}\cr
$\stackrel{}{{}_{q\to \infty}}$\cr
}} }
\sum\limits_{j_1, j_2, j_3=0}^{q}
C_{j_3 j_2 j_1}\zeta_{j_1}^{(i_1)}\zeta_{j_2}^{(i_2)}\zeta_{j_3}^{(i_3)}
\end{equation}

\vspace{4mm}
\noindent
converging in the mean-square sense is valid, where

$$
C_{j_3 j_2 j_1}=\int\limits_t^T\psi_3(t_3)\phi_{j_3}(t_3)
\int\limits_t^{t_3}\psi_2(t_2)\phi_{j_2}(t_2)
\int\limits_t^{t_2}\psi_1(t_1)\phi_{j_1}(t_1)dt_1dt_2dt_3;
$$

\vspace{4mm}
\noindent
another notations are the same as in Theorems {\rm 1, 2}.}

\vspace{2mm}

{\bf Theorem 6}\ \cite{2011-2}-\cite{xxx333}, 
\cite{2010-2}-\cite{2013}, \cite{271a},
\cite{arxiv-5}, \cite{arxiv-6}.\ 
{\it Suppose that 
$\{\phi_j(x)\}_{j=0}^{\infty}$ is a complete orthonormal system of 
Legendre polynomials or trigonometric functions in $L_2([t, T]).$
Then, for the iterated Stra\-to\-no\-vich stochastic integral
of multiplicity {\rm 4}

$$
{\int\limits_t^{*}}^T
{\int\limits_t^{*}}^{t_4}
{\int\limits_t^{*}}^{t_3}
{\int\limits_t^{*}}^{t_2}
d{\bf w}_{t_1}^{(i_1)}
d{\bf w}_{t_2}^{(i_2)}d{\bf w}_{t_3}^{(i_3)}d{\bf w}_{t_4}^{(i_4)}
$$

\vspace{3mm}
\noindent
the following 
expansion

$$
{\int\limits_t^{*}}^T
{\int\limits_t^{*}}^{t_4}
{\int\limits_t^{*}}^{t_3}
{\int\limits_t^{*}}^{t_2}
d{\bf w}_{t_1}^{(i_1)}
d{\bf w}_{t_2}^{(i_2)}d{\bf w}_{t_3}^{(i_3)}d{\bf w}_{t_4}^{(i_4)}=
\hbox{\vtop{\offinterlineskip\halign{
\hfil#\hfil\cr
{\rm l.i.m.}\cr
$\stackrel{}{{}_{q\to \infty}}$\cr
}} }
\sum\limits_{j_1, j_2, j_3, j_4=0}^{q}
C_{j_4 j_3 j_2 j_1}\zeta_{j_1}^{(i_1)}\zeta_{j_2}^{(i_2)}\zeta_{j_3}^{(i_3)}
\zeta_{j_4}^{(i_4)}
$$

\vspace{4mm}
\noindent
converging in the mean-square sense  
is valid, where $i_1, i_2, i_3, i_4=0, 1,\ldots,m,$

\vspace{1mm}
$$
C_{j_4 j_3 j_2 j_1}=\int\limits_t^T\phi_{j_4}(t_4)\int\limits_t^{t_4}
\phi_{j_3}(t_3)
\int\limits_t^{t_3}\phi_{j_2}(t_2)\int\limits_t^{t_2}\phi_{j_1}(t_1)
dt_1dt_2dt_3dt_4;
$$

\vspace{4mm}
\noindent
another notations are the same as in Theorems {\rm 1, 2}.}

\vspace{2mm}

Recently, a new approach to the expansion and mean-square 
approximation of iterated Stratonovich stochastic integrals has been obtained
\cite{2018a} (Sect.~2.10--2.16), \cite{arxiv-4} (Sect.~7--13), \cite{arxiv-5} (Sect.~13--19), 
\cite{arxiv-6} (Sect.~5--11), \cite{new-art-1-xxy}
(Sect.~4--9).
Let us formulate four theorems that were obtained using this approach.

\vspace{2mm}

{\bf Theorem 7}\ \cite{2018a}, \cite{arxiv-4}, \cite{arxiv-5}, 
\cite{arxiv-6}, \cite{new-art-1-xxy}.\
{\it Suppose 
that $\{\phi_j(x)\}_{j=0}^{\infty}$ is a complete orthonormal system of 
Legendre polynomials or trigonometric functions in the space $L_2([t, T]).$
Furthermore, let $\psi_1(\tau), \psi_2(\tau), \psi_3(\tau)$ are continuously dif\-ferentiable 
nonrandom functions on $[t, T].$ 
Then, for the 
iterated Stra\-to\-no\-vich stochastic integral of third multiplicity

$$
J^{*}[\psi^{(3)}]_{T,t}={\int\limits_t^{*}}^T\psi_3(t_3)
{\int\limits_t^{*}}^{t_3}\psi_2(t_2)
{\int\limits_t^{*}}^{t_2}\psi_1(t_1)
d{\bf w}_{t_1}^{(i_1)}
d{\bf w}_{t_2}^{(i_2)}d{\bf w}_{t_3}^{(i_3)}\ \ \ (i_1,i_2,i_3=0,1,\ldots,m)
$$

\vspace{4mm}
\noindent
the following 
relations

\vspace{-1mm}
\begin{equation}
\label{fin1}
J^{*}[\psi^{(3)}]_{T,t}
=\hbox{\vtop{\offinterlineskip\halign{
\hfil#\hfil\cr
{\rm l.i.m.}\cr
$\stackrel{}{{}_{p\to \infty}}$\cr
}} }
\sum\limits_{j_1, j_2, j_3=0}^{p}
C_{j_3 j_2 j_1}\zeta_{j_1}^{(i_1)}\zeta_{j_2}^{(i_2)}\zeta_{j_3}^{(i_3)},
\end{equation}

\vspace{3mm}
\begin{equation}
\label{fin2}
{\sf M}\left\{\left(
J^{*}[\psi^{(3)}]_{T,t}-
\sum\limits_{j_1, j_2, j_3=0}^{p}
C_{j_3 j_2 j_1}\zeta_{j_1}^{(i_1)}\zeta_{j_2}^{(i_2)}\zeta_{j_3}^{(i_3)}\right)^2\right\}
\le \frac{C}{p}
\end{equation}

\vspace{5mm}
\noindent
are fulfilled, where $i_1, i_2, i_3=0,1,\ldots,m$ in {\rm (\ref{fin1})} and 
$i_1, i_2, i_3=1,\ldots,m$ in {\rm (\ref{fin2}),}
constant $C$ is independent of $p,$

$$
C_{j_3 j_2 j_1}=\int\limits_t^T\psi_3(t_3)\phi_{j_3}(t_3)
\int\limits_t^{t_3}\psi_2(t_2)\phi_{j_2}(t_2)
\int\limits_t^{t_2}\psi_1(t_1)\phi_{j_1}(t_1)dt_1dt_2dt_3
$$

\vspace{4mm}
\noindent
and
$$
\zeta_{j}^{(i)}=
\int\limits_t^T \phi_{j}(\tau) d{\bf f}_{\tau}^{(i)}
$$ 

\vspace{2mm}
\noindent
are independent standard Gaussian random variables for various 
$i$ or $j$ {\rm (}in the case when $i\ne 0${\rm );} 
another notations are the same as in Theorems~{\rm 1, 2}.}

\vspace{2mm}

{\bf Theorem 8}\ \cite{2018a}, \cite{arxiv-4}, \cite{arxiv-5}, 
\cite{arxiv-6}, \cite{new-art-1-xxy}.\ {\it Let
$\{\phi_j(x)\}_{j=0}^{\infty}$ be a complete orthonormal system of 
Legendre polynomials or trigonometric functions in the space $L_2([t, T]).$
Furthermore, let $\psi_1(\tau), \ldots, \psi_4(\tau)$ be continuously dif\-ferentiable 
nonrandom functions on $[t, T].$ 
Then, for the 
iterated Stra\-to\-no\-vich stochastic integral of fourth multiplicity

\begin{equation}
\label{fin0}
J^{*}[\psi^{(4)}]_{T,t}={\int\limits_t^{*}}^T\psi_4(t_4)
{\int\limits_t^{*}}^{t_4}\psi_3(t_3)
{\int\limits_t^{*}}^{t_3}\psi_2(t_2)
{\int\limits_t^{*}}^{t_2}\psi_1(t_1)
d{\bf w}_{t_1}^{(i_1)}
d{\bf w}_{t_2}^{(i_2)}d{\bf w}_{t_3}^{(i_3)}d{\bf w}_{t_4}^{(i_4)}
\end{equation}

\vspace{4mm}
\noindent
the following 
relations

\begin{equation}
\label{fin3}
J^{*}[\psi^{(4)}]_{T,t}
=\hbox{\vtop{\offinterlineskip\halign{
\hfil#\hfil\cr
{\rm l.i.m.}\cr
$\stackrel{}{{}_{p\to \infty}}$\cr
}} }
\sum\limits_{j_1, j_2, j_3,j_4=0}^{p}
C_{j_4j_3 j_2 j_1}\zeta_{j_1}^{(i_1)}\zeta_{j_2}^{(i_2)}\zeta_{j_3}^{(i_3)}\zeta_{j_4}^{(i_4)},
\end{equation}

\vspace{3mm}

\begin{equation}
\label{fin4}
{\sf M}\left\{\left(
J^{*}[\psi^{(4)}]_{T,t}-
\sum\limits_{j_1, j_2, j_3, j_4=0}^{p}
C_{j_4 j_3 j_2 j_1}\zeta_{j_1}^{(i_1)}\zeta_{j_2}^{(i_2)}\zeta_{j_3}^{(i_3)}
\zeta_{j_4}^{(i_4)}
\right)^2\right\}
\le \frac{C}{p^{1-\varepsilon}}
\end{equation}

\vspace{5mm}
\noindent
are fulfilled, where $i_1, \ldots , i_4=0,1,\ldots,m$ in {\rm (\ref{fin0}),} {\rm (\ref{fin3})} 
and $i_1, \ldots, i_4=1,\ldots,m$ in {\rm (\ref{fin4}),}
constant $C$ does not depend on $p,$
$\varepsilon$ is an arbitrary
small positive real number 
for the case of complete orthonormal system of 
Legendre polynomials in the space $L_2([t, T])$
and $\varepsilon=0$ for the case of
complete orthonormal system of 
trigonometric functions in the space $L_2([t, T]),$

$$
C_{j_4 j_3 j_2 j_1}=
$$

$$
=
\int\limits_t^T\psi_4(t_4)\phi_{j_4}(t_4)
\int\limits_t^{t_4}\psi_3(t_3)\phi_{j_3}(t_3)
\int\limits_t^{t_3}\psi_2(t_2)\phi_{j_2}(t_2)
\int\limits_t^{t_2}\psi_1(t_1)\phi_{j_1}(t_1)dt_1dt_2dt_3dt_4;
$$

\vspace{4mm}
\noindent
another notations are the same as in Theorem~{\rm 7}.}

\vspace{2mm}

{\bf Theorem 9}\ \cite{2018a}, \cite{arxiv-4}, \cite{arxiv-5}, 
\cite{arxiv-6}, \cite{new-art-1-xxy}.\
{\it Assume 
that $\{\phi_j(x)\}_{j=0}^{\infty}$ is a complete orthonormal system of 
Legendre polynomials or trigonometric functions in the space $L_2([t, T])$
and $\psi_1(\tau), \ldots, \psi_5(\tau)$ are continuously dif\-ferentiable 
nonrandom functions on $[t, T].$ 
Then, for the 
iterated Stra\-to\-no\-vich stochastic integral of fifth multiplicity

\begin{equation}
\label{fin7}
J^{*}[\psi^{(5)}]_{T,t}={\int\limits_t^{*}}^T\psi_5(t_5)
\ldots
{\int\limits_t^{*}}^{t_2}\psi_1(t_1)
d{\bf w}_{t_1}^{(i_1)}
\ldots d{\bf w}_{t_5}^{(i_5)}
\end{equation}

\vspace{4mm}
\noindent
the following 
relations

\begin{equation}
\label{fin8}
J^{*}[\psi^{(5)}]_{T,t}
=\hbox{\vtop{\offinterlineskip\halign{
\hfil#\hfil\cr
{\rm l.i.m.}\cr
$\stackrel{}{{}_{p\to \infty}}$\cr
}} }
\sum\limits_{j_1,\ldots,j_5=0}^{p}
C_{j_5 \ldots j_1}\zeta_{j_1}^{(i_1)}\ldots \zeta_{j_5}^{(i_5)},
\end{equation}

\vspace{3mm}

\begin{equation}
\label{fin9}
{\sf M}\left\{\left(
J^{*}[\psi^{(5)}]_{T,t}-
\sum\limits_{j_1, \ldots, j_5=0}^{p}
C_{j_5 \ldots j_1}\zeta_{j_1}^{(i_1)}\ldots
\zeta_{j_5}^{(i_5)}
\right)^2\right\}
\le \frac{C}{p^{1-\varepsilon}}
\end{equation}

\vspace{5mm}
\noindent
are fulfilled, where $i_1, \ldots , i_5=0,1,\ldots,m$ in {\rm (\ref{fin7}),} {\rm (\ref{fin8})} 
and $i_1, \ldots, i_5=1,\ldots,m$ in {\rm (\ref{fin9}),}
constant $C$ is independent of $p,$
$\varepsilon$ is an arbitrary
small positive real number 
for the case of complete orthonormal system of 
Legendre polynomials in the space $L_2([t, T])$
and $\varepsilon=0$ for the case of
complete orthonormal system of 
trigonometric functions in the space $L_2([t, T]),$

$$
C_{j_5 \ldots j_1}=
\int\limits_t^T\psi_5(t_5)\phi_{j_5}(t_5)\ldots
\int\limits_t^{t_2}\psi_1(t_1)\phi_{j_1}(t_1)dt_1\ldots dt_5;
$$

\vspace{3mm}
\noindent
another notations are the same as in Theorems~{\rm 7, 8}.}

\vspace{2mm}

{\bf Theorem 10}\ \cite{2018a}, \cite{arxiv-4}, \cite{arxiv-5}, 
\cite{arxiv-6}, \cite{new-art-1-xxy}.\
{\it Suppose that 
$\{\phi_j(x)\}_{j=0}^{\infty}$ is a complete orthonormal system of 
Legendre polynomials or trigonometric functions in the space $L_2([t, T]).$
Then, for the 
iterated Stratonovich stochastic integral of sixth multiplicity

\begin{equation}
\label{after10001qu1}
J_{T,t}^{*(i_1\ldots i_6)}={\int\limits_t^{*}}^T
\ldots
{\int\limits_t^{*}}^{t_2}
d{\bf w}_{t_1}^{(i_1)}
\ldots d{\bf w}_{t_6}^{(i_6)}
\end{equation}

\vspace{3mm}
\noindent
the following 
expansion 

\vspace{-1mm}
$$
J_{T,t}^{*(i_1\ldots i_6)}
=\hbox{\vtop{\offinterlineskip\halign{
\hfil#\hfil\cr
{\rm l.i.m.}\cr
$\stackrel{}{{}_{p\to \infty}}$\cr
}} }
\sum\limits_{j_1, \ldots, j_6=0}^{p}
C_{j_6 \ldots j_1}\zeta_{j_1}^{(i_1)}\ldots
\zeta_{j_6}^{(i_6)}
$$

\vspace{4mm}
\noindent
that converges in the mean-square sense is valid, where
$i_1, \ldots, i_6=0, 1,\ldots,m,$

$$
C_{j_6 \ldots j_1}=
\int\limits_t^T\phi_{j_6}(t_6)\ldots
\int\limits_t^{t_2}\phi_{j_1}(t_1)dt_1\ldots dt_6;
$$

\vspace{3mm}
\noindent
another notations are the same as in Theorems~{\rm 7--9}.}

\vspace{5mm}

\section{Approximation of Iterated 
Stratonovich Stochastic Integrals Based on Multiple
Fourier--Legendre Series}

\vspace{5mm}

As was mentioned above, 
one of the main problems arising in the implementation of the 
numerical scheme (\ref{4.470}) is the joint
numerical modeling of the iterated Stratonovich stochastic integrals 
figuring in (\ref{4.470}). Let us consider efficient 
numerical modeling formulas for 
the iterated Stratonovich stochastic integrals based on Theorems 4--9.

Using Theorems 1, 2 ($k=1$), Theorems 4--9, and multiple Fourier--Legendre series, 
we obtain the following 
approximations of iterated Stratonovich stochastic 
integrals from (\ref{4.470}) \cite{2006}-\cite{arxiv-6}

\vspace{2mm}
\begin{equation}
\label{ccc1}
I_{(0)\tau_{p+1},\tau_p}^{*(i_1)}=\sqrt{\Delta}\zeta_0^{(i_1)},
\end{equation}

\vspace{3mm}

\begin{equation}
\label{ccc2}
I_{(1)\tau_{p+1},\tau_p}^{*(i_1)}=
-\frac{{\Delta}^{3/2}}{2}\left(\zeta_0^{(i_1)}+
\frac{1}{\sqrt{3}}\zeta_1^{(i_1)}\right),
\end{equation}

\vspace{3mm}

\begin{equation}
\label{ccc3}
{I}_{(2)\tau_{p+1},\tau_p}^{*(i_1)}=
\frac{\Delta^{5/2}}{3}\left(
\zeta_0^{(i_1)}+\frac{\sqrt{3}}{2}\zeta_1^{(i_1)}+
\frac{1}{2\sqrt{5}}\zeta_2^{(i_1)}\right),
\end{equation}

\vspace{4mm}
\begin{equation}
\label{ccc4}
I_{(00)\tau_{p+1},\tau_p}^{*(i_1 i_2)q}=
\frac{\Delta}{2}\left(\zeta_0^{(i_1)}\zeta_0^{(i_2)}+\sum_{i=1}^{q}
\frac{1}{\sqrt{4i^2-1}}\left(
\zeta_{i-1}^{(i_1)}\zeta_{i}^{(i_2)}-
\zeta_i^{(i_1)}\zeta_{i-1}^{(i_2)}\right)\right),
\end{equation}

\vspace{6mm}

$$
I_{(01)\tau_{p+1},\tau_p}^{*(i_1 i_2)q}=
-\frac{\Delta}{2}
I_{(00)\tau_{p+1},\tau_p}^{*(i_1 i_2)q}
-\frac{{\Delta}^2}{4}\Biggl(
\frac{1}{\sqrt{3}}\zeta_0^{(i_1)}\zeta_1^{(i_2)}+\Biggr.
$$

\vspace{3mm}
\begin{equation}
\label{ccc5}
+\Biggl.\sum_{i=0}^{q}\Biggl(
\frac{(i+2)\zeta_i^{(i_1)}\zeta_{i+2}^{(i_2)}
-(i+1)\zeta_{i+2}^{(i_1)}\zeta_{i}^{(i_2)}}
{\sqrt{(2i+1)(2i+5)}(2i+3)}-
\frac{\zeta_i^{(i_1)}\zeta_{i}^{(i_2)}}{(2i-1)(2i+3)}\Biggr)\Biggr),
\end{equation}

\vspace{6mm}

$$
I_{(10)\tau_{p+1},\tau_p}^{*(i_1 i_2)q}=
-\frac{\Delta}{2}I_{(00)\tau_{p+1},\tau_p}^{*(i_1 i_2)q}
-\frac{\Delta^2}{4}\Biggl(
\frac{1}{\sqrt{3}}\zeta_0^{(i_2)}\zeta_1^{(i_1)}+\Biggr.
$$

\vspace{3mm}
\begin{equation}
\label{ccc6}
+\Biggl.\sum_{i=0}^{q}\Biggl(
\frac{(i+1)\zeta_{i+2}^{(i_2)}\zeta_{i}^{(i_1)}
-(i+2)\zeta_{i}^{(i_2)}\zeta_{i+2}^{(i_1)}}
{\sqrt{(2i+1)(2i+5)}(2i+3)}+
\frac{\zeta_i^{(i_1)}\zeta_{i}^{(i_2)}}{(2i-1)(2i+3)}\Biggr)\Biggr)
\end{equation}

\vspace{5mm}
\noindent
or
\begin{equation}
\label{ccc50}
I_{(01)\tau_{p+1},\tau_p}^{*(i_1 i_2)q}
=
\sum_{j_1,j_2=0}^{q}
C_{j_2j_1}^{01}
\zeta_{j_1}^{(i_1)}\zeta_{j_2}^{(i_2)},
\end{equation}

\vspace{3mm}
\begin{equation}
\label{ccc51}
I_{(10)\tau_{p+1},\tau_p}^{*(i_1 i_2)q}
=
\sum_{j_1,j_2=0}^{p}
C_{j_2j_1}^{10}
\zeta_{j_1}^{(i_1)}\zeta_{j_2}^{(i_2)};
\end{equation}

\vspace{2mm}

\begin{equation}
\label{ccc7}
I_{(000)\tau_{p+1},\tau_p}^{*(i_1 i_2 i_3)q}
=\sum_{j_1, j_2, j_3=0}^{q}
C_{j_3 j_2 j_1}
\zeta_{j_1}^{(i_1)}\zeta_{j_2}^{(i_2)}\zeta_{j_3}^{(i_3)},
\end{equation}

\vspace{2mm}

\begin{equation}
\label{ccc8}
I_{(100)\tau_{p+1},\tau_p}^{*(i_1 i_2 i_3)q}
=\sum_{j_1, j_2, j_3=0}^{q}
C_{j_3 j_2 j_1}^{100}
\zeta_{j_1}^{(i_1)}\zeta_{j_2}^{(i_2)}\zeta_{j_3}^{(i_3)},
\end{equation}

\vspace{2mm}

\begin{equation}
\label{ccc9}
I_{(010)\tau_{p+1},\tau_p}^{*(i_1 i_2 i_3)q}
=\sum_{j_1, j_2, j_3=0}^{q}
C_{j_3 j_2 j_1}^{010}
\zeta_{j_1}^{(i_1)}\zeta_{j_2}^{(i_2)}\zeta_{j_3}^{(i_3)},
\end{equation}

\vspace{2mm}

\begin{equation}
\label{ccc10}
I_{(001)\tau_{p+1},\tau_p}^{*(i_1 i_2 i_3)q}
=\sum_{j_1, j_2, j_3=0}^{q}
C_{j_3 j_2 j_1}^{001}
\zeta_{j_1}^{(i_1)}\zeta_{j_2}^{(i_2)}\zeta_{j_3}^{(i_3)},
\end{equation}

\vspace{2mm}

\begin{equation}
\label{ccc11}
I_{(0000)\tau_{p+1},\tau_p}^{*(i_1 i_2 i_3 i_4)q}
=\sum_{j_1, j_2, j_3, j_4=0}^{q}
C_{j_4 j_3 j_2 j_1}
\zeta_{j_1}^{(i_1)}\zeta_{j_2}^{(i_2)}\zeta_{j_3}^{(i_3)}
\zeta_{j_4}^{(i_4)},
\end{equation}

\vspace{2mm}

\begin{equation}
\label{ccc12}
I_{(00000)\tau_{p+1},\tau_p}^{*(i_1 i_2 i_3 i_4 i_5)q}=
\sum\limits_{j_1, j_2, j_3, j_4, j_5=0}^{q}
C_{j_5j_4 j_3 j_2 j_1}\zeta_{j_1}^{(i_1)}\zeta_{j_2}^{(i_2)}\zeta_{j_3}^{(i_3)}
\zeta_{j_4}^{(i_4)}\zeta_{j_5}^{(i_5)},
\end{equation}

\vspace{7mm}

\noindent
where the Fourier--Legendre coefficients have the form

\vspace{1mm}
$$
C_{j_2j_1}^{01}
=
\int\limits_{\tau_p}^{\tau_{p+1}}(\tau_p-y)\phi_{j_3}(y)
\int\limits_{\tau_p}^{y}
\phi_{j_1}(x)dx dy =
\frac{\sqrt{(2j_1+1)(2j_2+1)}}{8}\Delta^{2}\bar
C_{j_2j_1}^{01},
$$

\vspace{3mm}

$$
C_{j_2j_1}^{10}
=\int\limits_{\tau_p}^{\tau_{p+1}}\phi_{j_3}(y)
\int\limits_{\tau_p}^{y}(\tau_p-x)
\phi_{j_1}(x)dx dy =
\frac{\sqrt{(2j_1+1)(2j_2+1)}}{8}\Delta^{2}\bar
C_{j_2j_1}^{10},
$$

\vspace{5mm}

$$
C_{j_3j_2j_1}=\int\limits_{\tau_p}^{\tau_{p+1}}\phi_{j_3}(z)
\int\limits_{\tau_p}^{z}\phi_{j_2}(y)
\int\limits_{\tau_p}^{y}
\phi_{j_1}(x)dx dy dz=
$$

\vspace{2mm}
\begin{equation}
\label{hhh1}
=
\frac{\sqrt{(2j_1+1)(2j_2+1)(2j_3+1)}}{8}\Delta^{3/2}\bar
C_{j_3j_2j_1},
\end{equation}

\vspace{6mm}

$$
C_{j_4j_3j_2j_1}=\int\limits_{\tau_p}^{\tau_{p+1}}\phi_{j_4}(u)
\int\limits_{\tau_p}^{u}\phi_{j_3}(z)
\int\limits_{\tau_p}^{z}\phi_{j_2}(y)
\int\limits_{\tau_p}^{y}
\phi_{j_1}(x)dx dy dz du=
$$

\vspace{2mm}
\begin{equation}
\label{hhh2}
=\frac{\sqrt{(2j_1+1)(2j_2+1)(2j_3+1)(2j_4+1)}}{16}\Delta^{2}\bar
C_{j_4j_3j_2j_1},
\end{equation}

\vspace{5mm}

$$
C_{j_3j_2j_1}^{001}=\int\limits_{\tau_p}^{\tau_{p+1}}(\tau_p-z)\phi_{j_3}(z)
\int\limits_{\tau_p}^{z}\phi_{j_2}(y)
\int\limits_{\tau_p}^{y}
\phi_{j_1}(x)dx dy dz=
$$

\vspace{2mm}
\begin{equation}
\label{hhh3}
=
\frac{\sqrt{(2j_1+1)(2j_2+1)(2j_3+1)}}{16}\Delta^{5/2}\bar
C_{j_3j_2j_1}^{001},
\end{equation}

\vspace{5mm}

$$
C_{j_3j_2j_1}^{010}=\int\limits_{\tau_p}^{\tau_{p+1}}\phi_{j_3}(z)
\int\limits_{\tau_p}^{z}(\tau_p-y)\phi_{j_2}(y)
\int\limits_{\tau_p}^{y}
\phi_{j_1}(x)dx dy dz=
$$

\vspace{2mm}
\begin{equation}
\label{hhh4}
=
\frac{\sqrt{(2j_1+1)(2j_2+1)(2j_3+1)}}{16}\Delta^{5/2}\bar
C_{j_3j_2j_1}^{010},
\end{equation}

\vspace{5mm}

$$
C_{j_3j_2j_1}^{100}=\int\limits_{\tau_p}^{\tau_{p+1}}\phi_{j_3}(z)
\int\limits_{\tau_p}^{z}\phi_{j_2}(y)
\int\limits_{\tau_p}^{y}
(\tau_p-x)\phi_{j_1}(x)dx dy dz=
$$

\vspace{2mm}
\begin{equation}
\label{hhh5}
=
\frac{\sqrt{(2j_1+1)(2j_2+1)(2j_3+1)}}{16}\Delta^{5/2}\bar
C_{j_3j_2j_1}^{100},
\end{equation}

\vspace{5mm}

$$
C_{j_5j_4 j_3 j_2 j_1}=
\int\limits_{\tau_p}^{\tau_{p+1}}\phi_{j_5}(v)
\int\limits_{\tau_p}^v\phi_{j_4}(u)
\int\limits_{\tau_p}^{u}
\phi_{j_3}(z)
\int\limits_{\tau_p}^{z}\phi_{j_2}(y)\int\limits_{\tau_p}^{y}\phi_{j_1}(x)
dxdydzdudv=
$$

\vspace{2mm}
\begin{equation}
\label{hhh6}
=\frac{\sqrt{(2j_1+1)(2j_2+1)(2j_3+1)(2j_4+1)(2j_5+1)}}{32}\Delta^{5/2}\bar
C_{j_5j_4 j_3 j_2 j_1},
\end{equation}

\vspace{5mm}
\noindent
where

\vspace{-2mm}

$$
\bar C_{j_2j_1}^{01}=-\int\limits_{-1}^{1}(1+y)P_{j_2}(y)
\int\limits_{-1}^{y}
P_{j_1}(x)dx dy,
$$

\vspace{1mm}

$$
\bar C_{j_2j_1}^{10}=-\int\limits_{-1}^{1}P_{j_2}(y)
\int\limits_{-1}^{y}
(1+x)P_{j_1}(x)dx dy,
$$

\vspace{1mm}

\begin{equation}
\label{jjj1}
\bar C_{j_3j_2j_1}=
\int\limits_{-1}^{1}P_{j_3}(z)
\int\limits_{-1}^{z}P_{j_2}(y)
\int\limits_{-1}^{y}
P_{j_1}(x)dx dy dz,
\end{equation}

\vspace{1mm}
\begin{equation}
\label{jjj2}
\bar C_{j_4j_3j_2j_1}=
\int\limits_{-1}^{1}P_{j_4}(u)
\int\limits_{-1}^{u}P_{j_3}(z)
\int\limits_{-1}^{z}P_{j_2}(y)
\int\limits_{-1}^{y}
P_{j_1}(x)dx dy dz,
\end{equation}

\vspace{1mm}

\begin{equation}
\label{jjj3}
\bar C_{j_3j_2j_1}^{100}=-
\int\limits_{-1}^{1}P_{j_3}(z)
\int\limits_{-1}^{z}P_{j_2}(y)
\int\limits_{-1}^{y}
P_{j_1}(x)(x+1)dx dy dz,
\end{equation}

\vspace{1mm}
\begin{equation}
\label{jjj4}
\bar C_{j_3j_2j_1}^{010}=-
\int\limits_{-1}^{1}P_{j_3}(z)
\int\limits_{-1}^{z}P_{j_2}(y)(y+1)
\int\limits_{-1}^{y}
P_{j_1}(x)dx dy dz,
\end{equation}

\vspace{1mm}

\begin{equation}
\label{jjj5}
\bar C_{j_3j_2j_1}^{001}=-
\int\limits_{-1}^{1}P_{j_3}(z)(z+1)
\int\limits_{-1}^{z}P_{j_2}(y)
\int\limits_{-1}^{y}
P_{j_1}(x)dx dy dz,
\end{equation}

\vspace{1mm}

\begin{equation}
\label{jjj6}
\bar C_{j_5j_4 j_3 j_2 j_1}=
\int\limits_{-1}^{1}P_{j_5}(v)
\int\limits_{-1}^{v}P_{j_4}(u)
\int\limits_{-1}^{u}P_{j_3}(z)
\int\limits_{-1}^{z}P_{j_2}(y)
\int\limits_{-1}^{y}
P_{j_1}(x)dx dy dz du dv,
\end{equation}

\vspace{4mm}
\noindent
where $P_i(x)$ $(i=0, 1, 2,\ldots)$ is the Legendre polynomial and

\vspace{2mm}

$$
\phi_i(x)=
\sqrt{\frac{2i+1}{\Delta}}P_i\left(\left(x-\tau_p-\frac{\Delta}{2}\right)
\frac{2}{\Delta}\right),\ \ \ i=0, 1, 2,\ldots 
$$

\vspace{5mm}

The Fourier--Legendre coefficients 

\vspace{1mm}
\begin{equation}
\label{sss1}
\bar C_{j_2 j_1}^{01},\ \bar C_{j_2 j_1}^{10},\
\bar C_{j_3 j_2 j_1},\ \bar C_{j_4 j_3 j_2 j_1},\ \bar C_{j_3 j_2 j_1}^{001},\ 
\bar C_{j_3 j_2 j_1}^{010},\ \bar C_{j_3 j_2 j_1}^{100},\
\bar C_{j_5 j_4 j_3 j_2 j_1}
\end{equation}

\vspace{4mm}
\noindent
can be calculated exactly before start of the numerical method (\ref{4.470}).
The above calculation can be done 
with Python, Derive or Maple.
In \cite{2006}, \cite{2017}-\cite{2013}, \cite{arxiv-3}
several tables with these coefficients can be found.
Moreover, in \cite{Kuz-Kuz}, \cite{Mikh-1} 
the database with 270,000 exactly calculated 
Fourier--Legendre coefficients (including (\ref{sss1})) 
was described.
This database was used in the software package,
which is written in the Python programming language
for the implementation of high-order strong
numerical schemes for Ito SDEs with non-commutative noise \cite{Kuz-Kuz}, \cite{Mikh-1}.
Note that the mentioned Fourier--Legendre coefficients
do not depend on the step of integration $\tau_{p+1}-\tau_p$ of the
numerical scheme,
which can be not a constant in a general case.

On the basis of the presented 
expansions (see (\ref{ccc1})--(\ref{ccc12})) of 
iterated Stratonovich stochastic integrals we 
can see that increasing of multiplicities of these integrals 
or degree indexes of their weight functions 
leads
to increasing 
of smallness orders with respect to $\Delta$ in the mean-square sense 
for iterated stochastic integrals.
This leads to a sharp decrease  
of member 
quantities (the numbers $q$)
in expansions of iterated Stratonovich stochastic 
integrals, which are required for achieving the acceptable accuracy
of approximation.
Generally speaking, the minimal values $q$ that guarantee the condition 
(\ref{ors})
for each approximation (\ref{ccc1})--(\ref{ccc12})
are various and abruptly decreasing with the growth of 
smallness orders with respect to $\Delta$ in the mean-square sense for
iterated stochastic integrals.

Consider in detail the question on calculation and estimation
of the mean-square approximation error for the iterated
Stratonovich stochastic integrals (\ref{str11}) (see \cite{2018a}, Chapter~5 for details). 

Let us consider the following iterated Ito stochastic integrals

\vspace{1mm}
$$
I_{(l_1\ldots \hspace{0.2mm}l_k)T,t}^{(i_1\ldots i_k)}
=
\int\limits_t^T
(t-t_k)^{l_k} \ldots \int\limits_t^{t_{2}}
(t-t_1)^{l_1} d{\bf f}_{t_1}^{(i_1)}\ldots
d{\bf f}_{t_k}^{(i_k)},
$$

\vspace{4mm}
\noindent
where $i_1,\ldots, i_k=1,\dots,m,$\ \  $l_1,\ldots,l_k=0, 1, 2,$\ \
$k=1, 2,\ldots, 5.$

According to the standard relations between iterated
Ito and Stratonovich stochastic integrals, we obtain w.~p.~1
(with probability 1)

\vspace{1mm}
$$
I_{(00)\tau_{p+1},\tau_p}^{(i_1 i_2)}=
I_{(00)\tau_{p+1},\tau_p}^{*(i_1 i_2)}-
\frac{1}{2}{\bf 1}_{\{i_1=i_2\}}\Delta,
$$

\vspace{2mm}
$$
I_{(10)\tau_{p+1},\tau_p}^{(i_1 i_2)}=
I_{(10)\tau_{p+1},\tau_p}^{*(i_1 i_2)}+
\frac{1}{4}{\bf 1}_{\{i_1=i_2\}}\Delta^2,
$$

\vspace{2mm}
$$
I_{(01)\tau_{p+1},\tau_p}^{(i_1 i_2)}=
I_{(01)\tau_{p+1},\tau_p}^{*(i_1 i_2)}+
\frac{1}{4}{\bf 1}_{\{i_1=i_2\}}\Delta^2.
$$

\vspace{5mm}

Moreover,
the mean-square approximation error
for the iterated Ito stochastic integral

\vspace{1mm}
$$
I_{(00)\tau_{p+1},\tau_p}^{(i_1 i_2)}\ \ \ (i_1\ne i_2)
$$

\vspace{5mm}
\noindent
equals to the 
mean-square approximation error
for the iterated Stratonovich stochastic integral (see \cite{2018a}, Sect.~5.1 
for details)

\vspace{1mm}
$$
I_{(00)\tau_{p+1},\tau_p}^{*(i_1 i_2)}\ \ \ (i_1\ne i_2).
$$

\vspace{7mm}

From Theorem 3 we obtain \cite{2006}-\cite{200a},
\cite{301a}-\cite{arxiv-6}

\vspace{2mm}
\begin{equation}
\label{1}
{\sf M}\Biggl\{\left(I_{(00)\tau_{p+1},\tau_p}^{(i_1 i_2)}-
I_{(00)\tau_{p+1},\tau_p}^{(i_1 i_2)q}
\right)^2\Biggr\}
=\frac{\Delta^2}{2}\Biggl(\frac{1}{2}-\sum_{i=1}^q
\frac{1}{4i^2-1}\Biggr)\ \ \ (i_1\ne i_2),
\end{equation}

\vspace{5mm}

$$
{\sf M}\Biggl\{\left(I_{(10)\tau_{p+1},\tau_p}^{(i_1 i_2)}-
I_{(10)\tau_{p+1},\tau_p}^{(i_1 i_2)
q}
\right)^2\Biggr\}=
{\sf M}\Biggl\{\left(I_{(01)\tau_{p+1},\tau_p}^{(i_1 i_2)}-
I_{(01)\tau_{p+1},\tau_p}^{(i_1 i_2)q}\right)^2\Biggr\}=
$$

\vspace{2mm}
$$
=\frac{\Delta^4}{16}\Biggl(\frac{5}{9}-
2\sum_{i=2}^{q}\frac{1}{4i^2-1}-
\sum_{i=1}^{q}
\frac{1}{(2i-1)^2(2i+3)^2}
-\Biggr.
$$

\vspace{2mm}
\begin{equation}
\label{2}
\Biggl.-
\sum_{i=0}^{q}\frac{(i+2)^2+(i+1)^2}{(2i+1)(2i+5)(2i+3)^2}
\Biggr)\ \ \ (i_1\ne i_2).
\end{equation}

\vspace{7mm}

The case $i_1=i_2$ is considered in \cite{2018a}, Sect.~5.1.

Let us estimate the mean-square approximation error
for the iterated Stratonovich stochastic integrals (\ref{str11})
of multiplicities $k\ge 3.$
From (\ref{1}) ($i_1\ne i_2$) we get

\vspace{2mm}
$$
{\sf M}\left\{\left(I_{(00)\tau_{p+1},\tau_p}^{*(i_1 i_2)}-
I_{(00)\tau_{p+1},\tau_p}^{*(i_1 i_2)q}
\right)^2\right\}=\frac{\Delta^2}{2}
\sum\limits_{i=q+1}^{\infty}\frac{1}{4i^2-1}\le 
$$

\vspace{3mm}
\begin{equation}
\label{teac}
\le \frac{\Delta^2}{2}\int\limits_{q}^{\infty}
\frac{1}{4x^2-1}dx
=-\frac{\Delta^2}{8}{\rm ln}\left|
1-\frac{2}{2q+1}\right|\le C_1\frac{\Delta^2}{q},
\end{equation}

\vspace{6mm}
\noindent
where constant $C_1$ does not depend on $\Delta$.

As was mentioned above,
the value $\Delta$ plays the role of integration step 
in the numerical procedures for Ito SDEs.
Then this value is a sufficiently small.
Keeping in mind this circumstance, it is easy to notice that there 
exists such a constant $C_2$ that

\vspace{1mm}
\begin{equation}
\label{teac3}
{\sf M}\left\{\left(I_{(l_1\ldots l_k)\tau_{p+1},\tau_p}^{*(i_1\ldots i_k)}-
I_{(l_1\ldots l_k)\tau_{p+1},\tau_p}^{*(i_1\ldots i_k)q}\right)^2\right\}
\le C_2 {\sf M}\left\{\left(I_{(00)\tau_{p+1},\tau_p}^{*(i_1 i_2)}-
I_{(00)\tau_{p+1},\tau_p}^{*(i_1 i_2)q}\right)^2\right\},
\end{equation}

\vspace{5mm}
\noindent
where $I_{(l_1\ldots l_k)\tau_{p+1},\tau_p}^{*(i_1\ldots i_k)q}$
is the approximation of the iterated Stratonovich stochastic integral 
(\ref{str11}) for $k\ge 3$.

From (\ref{teac}) and (\ref{teac3}) we finally have

\vspace{1mm}
\begin{equation}
\label{teac4}
{\sf M}\left\{\left(I_{(l_1\ldots l_k)\tau_{p+1},\tau_p}^{*(i_1\ldots i_k)}-
I_{(l_1\ldots l_k)\tau_{p+1},\tau_p}^{*(i_1\ldots i_k)q}\right)^2\right\}
\le K \frac{\Delta^2}{q},
\end{equation}

\vspace{5mm}
\noindent
where constant $K$ does not depend on $\Delta.$ 

The same idea can be found in \cite{KlPl2} 
for the case of trigonometric functions.
Note that, in contrast to the estimate (\ref{teac4}), 
the constant $C$ in Theorems 7--9 does not depend on $q.$

Essentially more information about numbers $q$ can be obtained 
by another approach. We have

\vspace{1mm}
$$
I_{(l_1\ldots l_k)\tau_{p+1},\tau_p}^{*(i_1\ldots i_k)}=
I_{(l_1\ldots l_k)\tau_{p+1},\tau_p}^{(i_1\ldots i_k)}\ \ \ \hbox{w.~p.~1}
$$

\vspace{5mm}
\noindent
for pairwise different $i_1,\ldots,i_k=1,\ldots,m$.

Then, for pairwise different $i_1,\ldots,i_5=1,\ldots,m$ 
from (\ref{qq1}) we obtain

\vspace{1mm}
$$
{\sf M}\left\{\left(
I_{(01)\tau_{p+1},\tau_p}^{*(i_1i_2)}-
I_{(01)\tau_{p+1},\tau_p}^{*(i_1i_2)q}\right)^2\right\}=
\frac{\Delta^{4}}{4}-\sum_{j_1,j_2=0}^{q}
\left(C_{j_2j_1}^{01}\right)^2,
$$

\vspace{3mm}
$$
{\sf M}\left\{\left(
I_{(10)\tau_{p+1},\tau_p}^{*(i_1i_2)}-
I_{(10)\tau_{p+1},\tau_p}^{*(i_1i_2)q}\right)^2\right\}=
\frac{\Delta^{4}}{12}-\sum_{j_1,j_2=0}^{q}
\left(C_{j_2j_1}^{10}\right)^2,
$$

\vspace{3mm}
$$
{\sf M}\left\{\left(
I_{(000)\tau_{p+1},\tau_p}^{*(i_1i_2 i_3)}-
I_{(000)\tau_{p+1},\tau_p}^{*(i_1i_2 i_3)q}\right)^2\right\}=
\frac{\Delta^{3}}{6}-\sum_{j_3,j_2,j_1=0}^{q}
C_{j_3j_2j_1}^2,
$$

\vspace{3mm}
$$
{\sf M}\left\{\left(
I_{(0000)\tau_{p+1},\tau_p}^{*(i_1i_2 i_3 i_4)}-
I_{(0000)\tau_{p+1},\tau_p}^{*(i_1i_2 i_3 i_4)q}\right)^2\right\}=
\frac{\Delta^{4}}{24}-\sum_{j_1,j_2,j_3,j_4=0}^{q}
C_{j_4j_3j_2j_1}^2,
$$

\vspace{3mm}
$$
{\sf M}\left\{\left(
I_{(100)\tau_{p+1},\tau_p}^{*(i_1i_2 i_3)}-
I_{(100)\tau_{p+1},\tau_p}^{*(i_1i_2 i_3)q}\right)^2\right\}=
\frac{\Delta^{5}}{60}-\sum_{j_1,j_2,j_3=0}^{q}
\left(C_{j_3j_2j_1}^{100}\right)^2,
$$

\vspace{3mm}
$$
{\sf M}\left\{\left(
I_{(010)\tau_{p+1},\tau_p}^{*(i_1i_2 i_3)}-
I_{(010)\tau_{p+1},\tau_p}^{*(i_1i_2 i_3)q}\right)^2\right\}=
\frac{\Delta^{5}}{20}-\sum_{j_1,j_2,j_3=0}^{q}
\left(C_{j_3j_2j_1}^{010}\right)^2,
$$

\vspace{3mm}
$$
{\sf M}\left\{\left(
I_{(001)\tau_{p+1},\tau_p}^{*(i_1i_2 i_3)}-
I_{(001)\tau_{p+1},\tau_p}^{*(i_1i_2 i_3)q}\right)^2\right\}=
\frac{\Delta^5}{10}-\sum_{j_1,j_2,j_3=0}^{q}
\left(C_{j_3j_2j_1}^{001}\right)^2,
$$

\vspace{3mm}
$$
{\sf M}\left\{\left(
I_{(00000)\tau_{p+1},\tau_p}^{*(i_1 i_2 i_3 i_4 i_5)}-
I_{(00000)\tau_{p+1},\tau_p}^{*(i_1 i_2 i_3 i_4 i_5)q}\right)^2\right\}=
\frac{\Delta^{5}}{120}-\sum_{j_1,j_2,j_3,j_4,j_5=0}^{q}
C_{j_5 i_4 i_3 i_2 j_1}^2.
$$

\vspace{7mm}

For example \cite{2006}-\cite{2013},

\vspace{2mm}

$$
{\sf M}\left\{\left(
I_{(000)\tau_{p+1},\tau_p}^{*(i_1i_2 i_3)}-
I_{(000)\tau_{p+1},\tau_p}^{*(i_1i_2 i_3)6}\right)^2\right\}=
\frac{\Delta^{3}}{6}-\sum_{j_3,j_2,j_1=0}^{6}
C_{j_3j_2j_1}^2
\approx
0.01956000\Delta^3,
$$

\vspace{3mm}
$$
{\sf M}\left\{\left(
I_{(100)\tau_{p+1},\tau_p}^{*(i_1i_2 i_3)}-
I_{(100)\tau_{p+1},\tau_p}^{*(i_1i_2 i_3)2}\right)^2\right\}=
\frac{\Delta^{5}}{60}-\sum_{j_1,j_2,j_3=0}^{2}
\left(C_{j_3j_2j_1}^{100}\right)^2
\approx
0.00815429\Delta^5,
$$

\vspace{3mm}
$$
{\sf M}\left\{\left(
I_{(010)\tau_{p+1},\tau_p}^{*(i_1i_2 i_3)}-
I_{(010)\tau_{p+1},\tau_p}^{*(i_1i_2 i_3)2}\right)^2\right\}=
\frac{\Delta^{5}}{20}-\sum_{j_1,j_2,j_3=0}^{2}
\left(C_{j_3j_2j_1}^{010}\right)^2
\approx
0.01739030\Delta^5,
$$

\vspace{3mm}
$$
{\sf M}\left\{\left(
I_{(001)\tau_{p+1},\tau_p}^{*(i_1i_2 i_3)}-
I_{(001)\tau_{p+1},\tau_p}^{*(i_1i_2 i_3)2}\right)^2\right\}=
\frac{\Delta^5}{10}-\sum_{j_1,j_2,j_3=0}^{2}
\left(C_{j_3j_2j_1}^{001}\right)^2
\approx 0.02528010\Delta^5,
$$

\vspace{3mm}
$$
{\sf M}\left\{\left(
I_{(0000)\tau_{p+1},\tau_p}^{*(i_1i_2i_3 i_4)}-
I_{(0000)\tau_{p+1},\tau_p}^{*(i_1i_2i_3 i_4)2}\right)^2\right\}=
\frac{\Delta^{4}}{24}-\sum_{j_1,j_2,j_3,j_4=0}^{2}
C_{j_4 j_3 j_2 j_1}^2\approx
0.02360840\Delta^4,
$$

\vspace{3mm}
$$
{\sf M}\left\{\left(
I_{(00000)\tau_{p+1},\tau_p}^{*(i_1i_2i_3i_4 i_5)}-
I_{(00000)\tau_{p+1},\tau_p}^{*(i_1i_2i_3i_4 i_5)1}\right)^2\right\}=
\frac{\Delta^5}{120}-\sum_{j_1,j_2,j_3,j_4,j_5=0}^{1}
C_{j_5 i_4 i_3 i_2 j_1}^2\approx
0.00759105\Delta^5.
$$

\end{document}